\documentclass{amsart}
\usepackage{amsmath,amssymb}
\usepackage{graphicx,calc,amscd}
\usepackage{epsfig,psfig}
\usepackage{float}
\usepackage{labelfig}
\usepackage{amsfonts}
\usepackage{latexsym}

\newtheorem{theorem}{Theorem}[section]
\newtheorem{lemma}[theorem]{Lemma}
\newtheorem{corollary}[theorem]{Corollary}
\newtheorem{proposition}[theorem]{Proposition}

\newtheorem{conjecture}{Conjecture}

\theoremstyle{definition}

\theoremstyle{remark}
\newtheorem{rem}{Remark}[section]

\numberwithin{equation}{section}

\newcommand{\ii}{\mathrm{int}}
\newcommand{\Z}{{\mathbb Z}}
\newcommand{\arcsinh}{\mathrm{arcsinh}}
\newcommand{\arccosh}{\mathrm{arccosh}}
\def\fullref{\ref}

\def\hyp{\Bbb H}

\def\reals{\Bbb R}
\def\RR{\Bbb R}

\def\twist{\bold T}

\def\deform{{\mathcal T}}
\def\model{{\mathcal T}^{\mbox{s}}_{1,1}}
\def\arc{{\mathcal A}}
\def\equal{{\mathcal E}}
\def\length{\ell}
\def\trace{\mbox{tr\,}}
\def\Neq{{\mathcal Neq}}
\def\complex{\Bbb C}
\def\sl{\mbox{SL}}

\def\psl{PSL(2,\RR)}

\def\NN{\Bbb N}
\def\cc{\gamma_{1}}

\begin{document}

\title[Hypersurfaces in Teichm\"uller Space]
{Multiplicities of simple closed geodesics and hypersurfaces in Teichm\"uller space}

\author{Greg McShane and Hugo Parlier}

\address{Laboratoire Emile Picard\\
Université Paul Sabatier\\
Toulouse\\
France}
\address{Section de Math\'ematiques\\
Universit\'e de Gen\`eve\\
Gen\`eve\\
Switzerland} \email{greg.mcshane@gmail.com}
\email{hugo.parlier@math.unige.ch}

\date{\today}

\maketitle

\begin{abstract}
Using geodesic length functions, we define a natural family of
real codimension 1 subvarieties of Teichm\"uller space, namely the
subsets where the lengths of two distinct simple closed geodesics
are of equal length. We investigate the point set topology of the
union of all such hypersurfaces using elementary methods. Finally,
this analysis is applied to investigate the nature of the Markoff
conjecture.
\end{abstract}

\section{Introduction}
We define the {\em simple length spectrum} of a Riemann surface of
genus $g$ with $n$ totally geodesic boundary components to be the
set of lengths of simple closed geodesics counted with
multiplicities. As the metric varies, the length spectrum changes.
We are interested in three questions.\\

Is there a surface for which all the multiplicities are 1?\\

How big is the set of such surfaces?\\

Is it possible to deform a surface such that the multiplicity
stays 1 for all simple geodesics?\\

The answer to the first question is of course positive, although
to answer it we shall show that the set of surfaces where the
spectrum does not have this property is Baire meagre. The answer
to the third question is clearly yes if we allow deformations in
the space of all Riemannian metrics, but we show that the answer
is no if we restrict ourselves to metrics of fixed constant
curvature.

The hypersurfaces we study are the non-empty subsets
$\equal(\alpha,\beta)$ of Teichm\"uller space where a pair of
distinct simple closed geodesics $\alpha,\beta$ have the same
length. When the intersection number $\ii(\alpha,\beta)$ is small,
the surfaces $\equal(\alpha,\beta)$ play an important part in the
theory of fundamental domains for the mapping class group in low
genus, see for instance \cite{gri98} and \cite{mas96}.\\

Our main theorem is the following:

\begin{theorem}\label{thm:noarc}
The set of surfaces with simple simple length spectrum is dense
and its complement is Baire meagre.

If $\arc$ is a path in Teichm\"uller space $\deform$ then there is
a surface on $\arc$ which has at least two distinct simple closed
geodesics of the same length.
\end{theorem}

Let $\equal$ denote the set of all surfaces with at least one pair
of simple closed geodesics of equal length; $\equal$ is the union
of a countable family of nowhere dense subsets, namely the sets
$\equal(\alpha,\beta)$ where $\alpha,\beta$ vary over all distinct
simple closed geodesics. The theorem asserts that $\equal$ is
dense and moreover that the complement contains no arcs and is
thus totally disconnected. Our next point of interest is the
topology of the sets $\equal(\alpha,\beta)$ and $\equal$.

\begin{theorem}\label{thm:B}
The sets $\equal(\alpha,\beta)$ are connected sub-manifolds of
Teichm\"uller space. The set $\equal$ is connected.
\end{theorem}

In the case of the Teichm\"uller space $\deform_{1,1}^\varepsilon$
of a one-holed torus with fixed boundary length, a careful study
of the asymptotic behaviour of the sets $\equal(\alpha,\beta)$
will allow us to prove the following.

\begin{theorem}\label{thm:C}
Let $\alpha$ and $\beta$ be a pair of distinct simple closed
geodesics on a one-holed torus. The set $\equal(\alpha,\beta)$ is
a simple path joining the points in the Thurston boundary of
$\deform_{1,1}^\varepsilon$ determined by the two unique simple
closed geodesics with equal number of intersection points with
both $\alpha$ and $\beta$.
\end{theorem}

We strongly suspect that the set of one-holed tori which contain
three simple closed geodesics of equal length is nowhere dense.

Though the Euclidean torus is not negatively curved, it is useful
to bear it in mind as a prototype. The reader is invited to check
that our theorems above hold for Euclidean tori. Let $\deform$ be
the Teichm\"uller space of Euclidean tori, which one identifies
with the upper half plane, or the $\tau$ plane, in the usual way.
The rational $p/q \in \RR$ corresponds to a simple curve of slope
$p/q$ on the square torus which we shall use as a reference
surface. The length of the curve on the surface corresponding to
the parameter $\tau$ is $|p\tau +q|$ and geometry of the set
$\equal$ is readily understood in terms of elementary number
theory. By computation, one sees that the set
$\equal(\alpha,\beta)$ is always a Poincar\'e geodesic, joining
either pairs of rationals or pairs of quadratic irrationals on the
real line. The map $z\mapsto -\overline{z} $ preserves the
rationals and fixes $[0,\infty]= \{ \Re z = 0 \}$. For any
rational $m/n$ the curves of slope $m/n$ and $-m/n$ have the same
length at $\tau\in [0,\infty]$. One maps $[0,\infty]$ onto any
geodesic $[a/b,a'/b']$ using $\psl$, which is transitive on pairs
of rationals and thus finds a pair of curves which are of equal
length on $[a/b,a'/b']$. The first sentence of Theorem
\ref{thm:noarc} follows since quadratic irrationals are countable
and the second because every pair of distinct points in the plane
is separated by a Poincar\'e geodesic with rational endpoints.
Unfortunately, in the more general situation we consider below, no
such characterization of the sets $\equal(\alpha,\beta)$ exists,
although we obtain an analogous result for one holed hyperbolic
tori with fixed boundary length in Theorem \ref{thm:C}.

It is also interesting to note that the simple length spectrum of
the square flat torus has unbounded multiplicity. The squares of
the lengths are sums of squares of integers. By elementary number
theory, for any sequence of integers $k_n$, each a product of $n$
distinct primes congruent to 1 mod 4, the number of ways of
writing $k_n$ as a sum as squares of coprime integers goes to
infinity with $n$. It is well known that a prime $p\neq 2$ can be
written as a sum of squares $x^2 + y^2$, $x,y\in \mathbb{N}$ if
and only if $p$ is congruent to $1$ modulo $4$. It is also well
known that there are infinitely many primes congruent to $1$
modulo $4$. Such a prime $p$ admits a factorization $p=
(x+iy)(x-iy)$ and $x+iy, x-iy$ are irreducible elements of the
ring of Gaussian integers $\mathbb{Z}[i]$. Choose $n$ distinct
primes $p_{k} \in \mathbb{N},\,1 \leq k \leq n$, let $a_{k}\in
\mathbb{Z}[i]$ such that $p_{k} = a_{k}\bar{a_{k}}$ and let $N$
denote their product. It is immediate from the above that $N$
factorizes over $\mathbb{Z}[i]$
$$N = (a_{1}\bar{a}_{1})(a_{2}\bar{a}_{2})\ldots (a_{n}\bar{a}_{n}).$$
Consider the set $R_{N}$ of Gaussian integers of the form
$c_{1}c_{2}\ldots c_{n}$ where $c_{k}\in \{a_{k},\bar{a}_{k}\}$.
Note that the modulus squared of each element of $R_{n}$ is $N$.
It is easy to check, using the fact that the Gaussian integers is
a unique factorization domain, that $R_{N}$ contains exactly
$2^{n-1}$ distinct elements. Note further that if
$c_{1}c_{2}\ldots c_{n} \in R_{N}$ and $c_{1}c_{2}\ldots c_{n}=
x+iy$ then $x,y$ are coprime integers, for otherwise there is a
prime $p$ that divides $x,y$ hence $x+iy$ but this is impossible
as the $c_{i}$ are irreducible and no two $c_{i}$ are complex
conjugate.

Similarly, for the spectrum of lengths of \textit{all} closed
geodesics of hyperbolic surfaces, the multiplicities are
unbounded. More precisely, Randol \cite{ra80} shows (using results
of Horowitz) that for any surface of constant curvature and any
$n>0$ there is a set of $n$ distinct primitive geodesics of the
same length. By the above, these curves necessarily have double
points.

Finally, by the work of Schmutz-Schaller \cite{sc98} and others on
the systole one has a lower bound on the multiplicity of the
simple length spectrum as the surface varies over Teichm\"uller
space. The existence (or not) of an upper bound for the
multiplicity of the simple length spectrum for surfaces of a given
signature is an open question. In particular, Schmutz-Schaller
conjectured the following in \cite{sc98} for once-punctured tori.

\begin{conjecture}\label{conj:schmutz} The multiplicity of the simple length
spectrum is bounded above by 6 for a once-punctured torus.
\end{conjecture}

In fact, this conjecture is a geometric generalization of an
equivalent version to the following well known conjecture, as will
be explained below.

\begin{conjecture}\label{conj:mark}
{\bf (Markoff conjecture)} The Markoff conjecture \cite{fr13} is a
conjecture in classical number theory. In its original form, it
concerns solutions in the positive integers to the cubic
$$x^2 +y^2 +z^2 - xyz =0.$$
It states simply that if there are solutions $x \geq y \geq z$ and
$x \geq y' \geq z'$ over the positive integers then $y=y'$ and
$z=z'$.
\end{conjecture}

At the time of writing, conjecture \ref{conj:mark} is known to be
correct for $x$ a prime power, see \cite{sc965}, \cite{butt98},
and more recently \cite{ser05}. Other results related to this
conjecture include \cite{ba96} and \cite{butt01}.

Over the past few decades, several authors (i.e. \cite{se85},
\cite{ha86}, \cite{sch03}, \cite{co55}, \cite{co71}, \cite{sa91})
have had much success in translating problems in number theory
into the language of surfaces and geodesics. Their efforts have in
some cases provided a greater clarity to the understanding of some
classical theorems, particularly in Diophantine approximation.

Let us outline the well-known correspondence between Markoff
triples $(x,y,z)$ (solutions to the above problem) and triples of
simple closed geodesics $(\alpha,\beta,\gamma)$ on a certain
once-punctured torus (a torus with a single cusp as boundary).

There is, up to isometry, exactly one hyperbolic once-punctured
torus, $\mathbf{M}$, with isometry group of order $12$. This
surface, commonly called the modular torus, is obtained as the
quotient of the hyperbolic plane by a certain subgroup, $G$ of
index 6 of the modular group $PSL(2,\Z)$ acting by linear
fractional transformations. In fact $G$ is a normal subgroup of
the modular group and the quotient group $SL(2,\Z) / G  \simeq C_2
\times S_3$ is isomorphic to the group of orientation preserving
isometries of the quotient surface. By Riemann-Hurwitz, 6 is the
maximum number of orientation preserving isometries of a
once-punctured torus and, by a similar argument, 12 is the maximum
number of isometries of a once-punctured torus.

If $(\alpha,\beta,\gamma)$ is a triple of simple geodesics on
$\mathbf{M}$ which pairwise intersect in exactly one point (i.e.
the three sets $\alpha \cap
\beta,\beta\cap\gamma,\gamma\cap\alpha$ each consist of a single
point) then
$$\left(2\cosh\frac{\length(\alpha)}{2},\\
 2\cosh\frac{\length(\beta)}{2},\\
 2\cosh\frac{\length(\gamma)}{2}\right)$$
is a solution to the cubic above. This is because one can find
matrices $A,B,C=AB$ such that any pair generate $G$ and so that
$\trace A = 2\cosh\frac{\length(\alpha)}{2}, \trace B =
2\cosh\frac{\length(\beta)}{2}, \trace C
=2\cosh\frac{\length(\gamma)}{2}$; the commutator of $A,B$ is {\it
parabolic}, i.e., has trace $-2$, and, using the familiar trace
identity in $SL(2,\complex)$ for the commutator \cite{beabook},
one sees that these three numbers satisfy the cubic.


There is an equivalent conjecture to the Markoff conjecture which
concerns simple closed geodesics on $\mathbf{M}$.

\begin{conjecture}\label{conj:mod}
The modular torus $\mathbf{M}$ has the following property: if
$\alpha,\beta$ are a pair of simple closed geodesics of the same
length, then there is an isometry of $\mathbf{M}$ taking one to
the other. We shall refer to this as the {\it Markoff isometry
property}.
\end{conjecture}

The equivalence of these conjectures can be shown as follows. We
begin by noting that given 3 real numbers $x,y,z$ satisfying the
cubic above there exist 3 matrices $A, B, AB=C \in SL(2,\reals)$
such that
$$\trace A = x, \trace B = y, \trace C=z,$$
and so that the commutator of $A,B$ is parabolic. Further, these
matrices are {\it unique up to conjugation} in $SL(2,\reals)$
\cite{go03}.

Now let $\alpha,\alpha'$ be a pair of simple geodesics on $M$
which have the same length. Since $\alpha$ (resp. $\alpha'$) is
simple there exist $A,B,C=AB \in G$ (resp. $A',B',C'=A'B' \in G$)
such that $A$ (resp. $A'$) covers $\alpha$ (resp. $\alpha'$), any
two of $A,B,C$ (resp. $A',B',C'$) generate $G$ and so that the
trace of the commutator of $A,B$ (resp. $A',B'$) is $-2$. From the
trace identity for this commutator, one sees that the quantities
$\trace A,\trace B,\trace C$ and $\trace A',\trace B',\trace C'$
satisfy the cubic. Since $\alpha$ and $\alpha'$ have the same
length $\trace A = \trace A'$ so assuming the Markoff conjecture
in its original form we must also have
$$\trace B = \trace B', \trace C = \trace C'.$$

By the uniqueness of $A,B,C$ in $SL(2,\reals)$ there exists a
matrix $T \in SL(2,\reals)$ so that $TAT^{-1}=A', TBT^{-1}=B',
TAT^{-1}=C'$. Thus conjugation by $T$ leaves $G$ invariant and so
this induces an isometry of $\hyp \backslash G$; since $T$
conjugates $A$ (which covers $\alpha$) to $A'$ (which covers
$\alpha'$) this isometry takes $\alpha$ to $\alpha'$ as
required.\\

Now $\mathbf{M}$ is the the only once-punctured torus which has
isometry group $C_2 \times S_3$, but not the only one-{\it holed}
torus. (When the boundary of torus is a cusp {\it or} the boundary
of the Nielsen core is a simple closed geodesic, we shall call it
a one-holed torus.) Choose a real number $t\geq 3$; there are
matrices, unique up to conjugacy in $SL(2,\reals)$, $A_t,B_t,
A_tB_t=C_t \in SL(2,\reals)$, each with trace $t$ and such that
the trace of the commutator of $A_t,B_t$ is $ 3t^2 -t^3 - 2$. The
group generated by $A_t,B_t$ is free, discrete and, when we think
of it as acting on the upper half space by linear fractional
transformations, the quotient surface, $M_t$, is homeomorphic to a
torus with a point removed. When $t=3$ the trace of the commutator
is $-2$ and the surface is the once-punctured torus $\mathbf{M}$.
When $t>3$ the commutator is {\it hyperbolic}, the {\it Nielsen
core} of the quotient is a proper subset and is bounded by a
simple closed geodesic which is covered by the commutator. For any
given $t\ge3$, from the uniqueness of the matrices,
$A_{t},B_{t},C_{t}$, it follows that the quotient surface has the
isometry group specified.

In view of this observation, one is led to consider the following
plausible generalization of the Markoff isometry property
conjecture:

\begin{conjecture}\label{conj:generalization}
Let $\model$ denote the set of hyperbolic one-holed tori which
have isometry group isomorphic to $C_2 \times S_3$. 
If $M$ is such a surface, then it enjoys the Markoff isometry
property; that is if $\alpha,\beta$ are a pair of simple closed
geodesics of the same length on $M$ then there is an isometry of
this surface taking one to the other.
\end{conjecture}

Notice that this conjecture is implied by a generalization to
one-{\it{holed}} tori of Schmutz-Schaller's conjecture mentioned
above. By applying the analysis used to prove Theorem
\fullref{thm:noarc}, we shall see that this conjecture is in fact
false, and thus a generalized version of Schmutz-Schaller's
conjecture is also false. More precisely, we show in Theorem
\ref{thm:dense} that the subset of surfaces $M$ which fail to have
the Markoff isometry property is dense in $\model$. This means
that the original Markoff conjecture is (probably) a conjecture in
pure number theory and not tractable by
hyperbolic geometry arguments.\\

\noindent{\bf Organization.}\\

We begin by reviewing some basic facts concerning curves and
geodesics from the theory of surfaces. Two preliminary results we
readily use, Theorems \ref{thm:mono1} and \ref{thm:mono2}, do not
seem to be known in the form we require, and because the proofs
are different in nature from the rest of the article, we defer
them to Appendix. There then follows a discussion of the {\it Dehn
twist} homeomorphism in relation to curves (lemma
\fullref{lem:twists}). Using Dehn twists, we generate special
sequences of geodesics and using these establish proposition
\fullref{prop:ratios}, which is fundamental to our argument. This
result, together with Theorem \fullref{thm:mono2}, implies Theorem
\fullref{thm:noarc}. We then discuss why we believe that a careful
modification of this argument provides evidence that the Markoff
conjecture cannot be proved using techniques from Teichm\"uller
theory. In the final section, we begin by showing that the
hypersurfaces we study are connected sub-manifolds and that their
union is also connected. Finally, a characterization of the
hypersurfaces is given in the particular case of one-holed tori
(with fixed boundary length), which is quite similar to the case
of flat tori mentioned above.\\

\noindent{\bf Acknowledgments.}\\

The first author would like to thank John Parker and Caroline
Series for introducing him to the Markoff conjecture, Andy Haas
for encouraging him to study the problem and David Epstein for
guidance and helpful remarks. The second author would like to
thank Benoit Bertrand, Thomas Gauglhofer and Alexander Lyatch for
precious advice. Both authors would like to thank a variety of
people they have bored with their unsuccessful geometric
constructions while trying to solve the Markoff conjecture,
including G\'erard Maze and Andr\'e Rocha.

\section{Preliminaries}

We review some elementary definitions and facts from the theory of
surfaces; much of this is available in a more detailed treatment
in either \cite{abbook}, \cite{beabook} or \cite{bubook}.
Throughout $M$ will denote a surface with constant curvature $-1$
and we shall insist that $M$ is {\it complete} with respect to
this metric (although we will only be concerned by what happens
inside the {\it Nielsen core} of the surface). This means that $M$
is {\it locally modeled} on the hyperbolic plane $\hyp^2$ and
there is a natural covering map $\pi:\hyp \rightarrow M$. By
$\deform$ we mean the Teichm\"uller space of $M$, meaning the
space of marked complete hyperbolic structures on $M$. The
signature of a surface $M$ will be denoted $(g,n)$ where $M$ is
homeomorphic to a surface of genus $g$ with $n$ simple closed
boundary curves. If deemed necessary, we will denote
$\deform_{g,n}$ the Teichm\"uller space of surfaces of signature
$(g,n)$. It is often useful to think of the Fenchel-Nielsen
coordinates of Teichm\"uller space, namely the lengths of a pants
decomposition of $M$ and a set of twist parameters which set how
the {\it interior} simple closed geodesics (defined below) are
pasted together. We insist on the fact the lengths of the boundary
curves (in fact the boundary curves of the Nielsen core) are
allowed to vary.\\

Next one recalls some facts concerning curves on surfaces (see
\cite{bubook} or \cite{cabl88} for details). Firstly, a {\it
simple curve} is a curve which has no self intersections. A curve
is {\it essential} if it bounds neither a disc nor a punctured
disc (annulus). For each free homotopy class which contains an
essential simple loop, there is a unique geodesic representative.
Further if such a free homotopy class contains a simple curve then
this geodesic representative is also simple. Simple closed
geodesics are called {\it separating} if they separate the surface
into two connected components and non-separating otherwise.
Furthermore, a simple closed geodesic is called {\it interior} if
it is entirely contained in the Nielsen core of the surface.

There is a natural function, $\length:\deform \times
\mbox{essential homotopy classes} \rightarrow \reals^+$, which
takes the pair $M,[\alpha]$ to the length $\ell_M(\alpha)$ of the
geodesic in the homotopy class $[\alpha]$ (measured in the
Riemannian metric on $M$). It is an abuse, though common in the
literature, to refer merely to {\it the length of the geodesic
$\alpha$} (rather than, more properly, the length of the geodesic
in the appropriate homotopy class).

The set of simple closed geodesics on the surface $M$ is more than
just an interesting curiosity. It was discovered by early
investigators (Fricke et al) that the lengths of a carefully
chosen finite subset of such curves could be used as a local
coordinate system for the space of surfaces; these are often
referred to as the {\it moduli} of the space. One way of seeing
this is through the following theorem, which we shall not use
explicitly although one could say that it contains the intuitive
idea of how we shall proceed.

\begin{theorem}\label{thm:mono1} Let $\deform$ be the Teichm\"uller
space of given signature. There is a fixed finite set of simple
closed geodesics $\gamma_1,\ldots,\gamma_n$ such that the map
$\varphi:$
$M\mapsto(\length_M(\gamma_1),\ldots,\length_M(\gamma_n))$ is
projectively injective on $\deform$.
\end{theorem}

(A map $f: X \rightarrow V$, $V$ a real vector space is {\it
projectively injective} if and only if $f(x) = t f(y),$ for some
$t \in
\reals$, implies $x=y$.)\\

It's important to note that in the case where $\deform$ is the
Teichm\"uller space of a surface with boundary, the geodesics are
allowed to be boundary geodesics. The nature of our investigation
requires us to study {\it interior} simple closed geodesics, so we
need a generalization of Theorem \ref{thm:mono1} where the set of
simple closed geodesics $\gamma_1,\ldots,\gamma_n$ are all
interior. While concocting a theorem suitable for our needs, we
discovered that the generalization to surfaces with boundary,
under the assumption that the geodesics are all interior, is
false. What is true however is the following theorem.

\begin{theorem}\label{thm:mono2} Let $\deform$ be the Teichm\"uller
space of given signature.

Then, there is a set $\Gamma:=\{\gamma_1,\hdots,\gamma_n\}$ of
{\it interior} simple closed curves such that for any given point
in $M\in\deform$, there are a finite many $M'\in\deform$ such that
$\lambda (\ell_{M}(\gamma_1),...,\ell_{M}(\gamma_n)) =
(\ell_{M'}(\gamma_1),...,\ell_{M'}(\gamma_n))$ for some $\lambda
\in\RR$.
\end{theorem}

We give proofs of both Theorem \ref{thm:mono1} and Theorem
\ref{thm:mono2}, but because these are different in nature to the
rest of the article, their proofs are deferred to Appendix. In the
Appendix, examples are given to highlight the differences between
Theorem \ref{thm:mono1} and Theorem \ref{thm:mono2}. Note that
Theorem \ref{thm:mono1} is probably known to specialists, and for
closed surfaces is implicit in \cite{th98}.

Of particular use to us, is the following corollary of Theorem
\ref{thm:mono2}.

\begin{corollary}\label{cor:ratio} Let $\deform$ be the Teichm\"uller space of given
signature, and $\Gamma$ be the set of curves of Theorem
\ref{thm:mono2}. For any $M\in
\deform$, there are finitely many $M'\in
\deform$ such that $\frac{\ell_{M'}(\alpha)}{\ell_{M'}(\beta)}
=\frac{\ell_{M}(\alpha)}{\ell_{M}(\beta)}$ for all $\alpha,\beta
\in \Gamma$.
\end{corollary}

If we consider only closed surfaces then the following argument
proves a slightly stronger result than is stated in our corollary
. Let $A,B$ be points in a Teichm\"uller space of a closed
hyperbolic surface; we now give a sketch of the existence of
simple closed curves $\alpha,\beta$
 such that
 $${\ell_A(\alpha)\over \ell_A(\beta)} \neq {\ell_B(\alpha)\over \ell_B(\beta)}.$$
 We shall suppose that there exists a simple closed curve $\gamma$
such that $A,B$ lie are in the same orbit for Fenchel-Nielsen
twist along $\gamma$; by replacing $\gamma$ by a measured
lamination and using Thurston's earthquake theorem one obtains the
general case. Suppose further that $\gamma$ is non separating then
there is a simple closed geodesic $\gamma'$ that meets $\gamma$ in
exactly one point. In fact, if $D_{\gamma}^n(\gamma')$ denotes the
simple geodesic freely homotopic to the image of $\gamma'$ by the
(right) Dehn twist along $\gamma$ iterated $n$-times, then
$D_{\gamma}^n(\gamma')$ meets $\gamma$ in exactly one point too.
Let $\theta_{n}$ denote the signed angle at
$D_{\gamma}^n(\gamma')\cap \gamma$, then it is not hard to show,
using hyperbolic trigonometry, that $\theta_{n}$ is monotone in
$n$, and furthermore that for any $\epsilon>0$ there exists $N$
such that
$$\cos(\theta_{-N})< -1+ \epsilon,\,\,\cos (\theta_{N}) >1- \epsilon.$$
Take $\alpha =D_{\gamma}^N(\gamma')$ and $\beta
=D_{\gamma}^N(\gamma')$. Following Kerkhoff \cite{ke83}, we see
that as we move from $A$ to $B$ along the Fenchel-Nielsen twist
orbit,
 $\ell(\alpha)$ is monotone decreasing
whilst $\ell_{\beta}$ is monotone increasing and
 $${\ell_A(\alpha)\over \ell_A(\beta)} < {\ell_B(\alpha)\over \ell_B(\beta)}.$$

Finally, with respect to the differential structure of
Teichm\"uller space, we have the following:

\begin{proposition}\label{prop:analytic}
For each closed geodesic $\alpha$, the function $M\mapsto
\length_M(\alpha)$ is an analytic function.
\end{proposition}

Probably the best way to prove this proposition is by means of
{\it Fricke trace calculus}. Briefly, one thinks of the surface as
being the quotient of $\hyp^2$ by a discrete group of isometries
$\Gamma$. To each closed geodesic one associates a conjugacy class
of elements in $\Gamma$ which {\it cover} it. If $A\in\Gamma$
covers $\alpha$ then $\trace A = 2 \cosh\frac{\length
(\alpha)}{2}$. Since $\Gamma$ is a subgroup of $\sl(2,\reals)$ one
can use the {\it trace relations} for this group to generate
relations between the lengths of geodesics on the surface. The
analytic nature of these relations is the key to showing the
proposition above (see for example \cite{abbook}).

\section{Dehn Twists and length ratios}\label{sec:ratios}

The proof of our main theorem \ref{thm:noarc} is based on a trick
involving the topological manipulation of curves. In particular,
we shall generate a sequence of simple closed curves with a
certain prescribed property which is outlined below (basically
their lengths should satisfy part 2 of lemma \ref{lem:twists}).

A simple topological argument, pasting together curves,
establishes the next proposition (left to the reader as an
exercise.)

\begin{lemma}\label{lem:twice}
Let $M$ be a surface and $\alpha$ a simple closed geodesic on $M$.
If $\alpha$ is not a boundary geodesic then there is a simple
closed geodesic $\alpha'$ which meets $\alpha$ in exactly 2
points.
\end{lemma}

\begin{rem}
We work with pairs of curves which meet in two points because a
simple closed geodesic which divides the surface into 2 components
meets any other geodesic loop in at least 2 points. A surface for
which every non peripheral simple closed geodesic is separating is
called a \textit{planar surface}. Topologically, planar surfaces
are $n$-holed spheres. The one-holed torus is not a planar
surface, and in fact it is the only topological type of surface
where every simple closed geodesic is non separating. By choosing
a pair of curves as in the lemma above, we are able to treat the
case of planar surfaces without much extra work.
\end{rem}

A Dehn twist $\twist_\alpha:M \mapsto M$ around a simple closed
curve $\alpha$ is a homeomorphism (defined up to isotopy) of the
surface $M$ to itself which is the identity on the complement of a
regular neighborhood of $\alpha$. The Dehn twist takes a geodesic
$\gamma$ meeting $\alpha$ to a (possibly non-geodesic) curve.
However, since this curve is also essential, we can
\textit{straighten it}, i.e., take the unique closed geodesic in
the homotopy class of $\twist_\alpha(\gamma)$. Since the twist is
a homeomorphism, the resulting geodesic is simple if and only if
the original geodesic $\gamma$ was. By convention, we say that a
Dehn twist takes simple closed geodesics to simple closed
geodesics.

This latter property of the Dehn twist round $\alpha$ allows us to
construct a sequence of geodesics with an interesting sequence of
lengths. We construct the sequence inductively as follows: let
$\alpha_0$ be any geodesic which meets $\alpha$ in exactly two
points as in lemma \ref{lem:twice}. Then let $\alpha_n$ be the
right Dehn twist of $\alpha_{n-1}$ round $\alpha$. This gives us a
well defined sequence of curves $\lbrace \alpha_n
\rbrace_{n=0}^{\infty}$, each of which is closed and simple.

\begin{lemma}\label{lem:twists}
With $\alpha,\lbrace \alpha_k \rbrace_{k=0}^{\infty}$ as above:
\begin{enumerate}
\item
For each $k=1,\ldots \infty$ the curve $\alpha_k$ meets $\alpha$
in 2 points and each $\alpha_k$ is simple.
\item
For all $k\in \NN$
$$2k\length(\alpha) -\length(\alpha_0) < \length (\alpha_k) \leq 2k \length(\alpha) + \length(\alpha_{0}).$$
\end{enumerate}
\end{lemma}

The first claim is immediate since homeomorphism takes simple
curves to simple curves; the second is essentially a consequence
of the triangle inequality. We note that, by using Riemannian
comparison theorems, this result is true even when $M$ does not
have {\it constant curvature} but merely {\it pinched negative
curvature}.

Fix a surface $M$ and let $\alpha,\beta$ be distinct simple closed
geodesics. Our aim is to calculate the ratio $\length(\alpha) /
\length(\beta)$ from the asymptotic formula in lemma
\ref{lem:twists}. For $\alpha$, resp. $\beta$, choose
$\alpha_{0}$, resp. $\beta_{0}$, as in lemma \ref{lem:twice}.
Define $\alpha_j := \twist_{\alpha}^j(\alpha_{0})$, resp. $\beta_j
:= \twist_{\beta}^j(\beta_{0})$). Set $B_{i}: = \lbrace \beta_{k}
: \length(\beta_k) \leq \length(\alpha_i)\rbrace.$

\begin{proposition}\label{prop:ratios}
With the notation above:
$$ {\sharp B_{i}\over i }\longrightarrow {\length(\alpha)\over \length(\beta)}.$$
\end{proposition}

\begin{proof}
By lemma \ref{lem:twists} we have
$$\sharp B_i \leq \sharp \{k : 2k \ell(\beta) - \ell(\beta_{0}) \leq 2i \ell(\alpha) + \ell(\alpha_{0}) \},$$
and
$$\sharp B_i \geq \sharp \{k : 2k \ell(\beta) + \ell(\beta_{0}) \leq 2i \ell(\alpha)  -  \ell(\alpha_{0}) \}.$$
It follows that
$$i{\ell(\alpha)\over \ell(\beta)} - {\ell(\beta_{0})\over 2 \ell(\beta)}
- {\ell(\alpha_{0})  \over 2\ell(\beta) }\leq \sharp B_i \leq
i{\ell(\alpha)\over \ell(\beta)} + {\ell(\alpha_{0})\over 2
\ell(\beta)}+ {\ell(\beta_{0})  \over 2\ell(\beta) } .$$ The
statement of the propostion is immediate.
\end{proof}

For the sequel, we shall denote $N(M)_i:=\sharp B_i$ for a given
choice of $M$.\\

\section{Surfaces with a pair of geodesics of equal length}

For completeness we give a brief account of the nature of the sets
$\equal(\alpha,\beta)$ in terms of elementary pointset topology.
The results presented in this section seem to be well known but do
not appear in the literature.

It should be possible to develop a very precise theory of these
sets in terms of algebraic functions, via Fricke trace calculus,
however, we will not consider that approach further here.

The next result exploits the fact, by proposition
\ref{prop:analytic}, that the function $M\mapsto
\length_M(\alpha)-\length_M(\beta)$ is analytic.

\begin{lemma}\label{lem:list}
Let $\alpha,\beta$ be a pair of simple closed geodesics. Then
\begin{enumerate}
\item
$\equal(\alpha,\beta)$ is non empty,
\item
$\equal(\alpha,\beta)$ is a closed subset with no interior, i.e.,
its complement is open dense in $\deform$,
\item
Let $\equal = \cup \equal(\alpha,\beta)$ where the union is taken
over all pairs of simple closed geodesics and set $\Neq = \deform
\backslash \equal$, then $\Neq$ is dense.
\end{enumerate}
\end{lemma}

\begin{proof} It is easy to show that, for each pair of simple closed
geodesics $\alpha \neq \beta$, $\equal(\alpha,\beta)$ is closed.
This follows when one notes that $\equal(\alpha,\beta)$ is the
zero set for the continuous function $M\mapsto
f_{\alpha,\beta}(M):=\length_M(\alpha)-\length_M(\beta)$.

Next one establishes that $\equal(\alpha,\beta)$ has no interior,
i.e., contains no open sets. Recall that if a real analytic
function is constant on an open set in its domain then in fact it
is constant on the entire component which contains this set. We
shall apply this to our function $f_{\alpha,\beta}$ to establish
the claim. Suppose that $\alpha$ and $\beta$ are simple geodesics
and they have the same length on some open subset of $\deform$ (so
that the function above is zero on this open set). Since $\deform$
is connected, the function is zero over the whole space, i.e., the
corresponding geodesics must have the same length over the whole
space. If $\alpha$ and $\beta$ are distinct then it is easy to
find a sequence of deformations of the surface such that the
length of $\alpha$ remains bounded while the length of $\beta$
tends to infinity (for example, if $\alpha$ meets $\beta$ the
right Dehn twist round $\alpha$ provides such a sequence). Thus
$\alpha=\beta$.

We show that $\equal(\alpha,\beta)$ is non empty by applying the
intermediate value theorem to our function $f_{\alpha,\beta}$ on
Teichm\"uller space, which is connected. The claim is true if
there are surfaces $X$ and $Y$ such that on $X$ $\alpha$ is very
short and $\beta$ is very long, and on $Y$ $\beta$ is very short
and $\alpha$ is very long. There are two cases to consider
depending on whether $\alpha$ and $\beta$ intersect or not.
Suppose that $\alpha$ and $\beta$ do not meet; then there is a
{\it pants decomposition} of $M$ so that both $\alpha$ and $\beta$
are curves in the decomposition. This enables one to find the
surfaces $X,Y$ quite easily, since one can specify the lengths of
the curves in a pants decomposition independently. On the other
hand, if $\alpha,\beta$ meet, then by the collar theorem
\cite{bubook}, for any $M\in \deform$, their lengths satisfy the
inequalities
\begin{eqnarray*}
\ell_M(\alpha) & \geq & 2
\arcsinh\left(\frac{1}{\sinh\frac{\ell_M(\beta)}{2}}\right),\\
\ell_M(\beta)  &\geq& 2
\arcsinh\left(\frac{1}{\sinh\frac{\ell_M(\alpha)}{2}}\right).
\end{eqnarray*}

Now by choosing a surface $X$, resp. $Y$, where the length of
$\alpha$ is indeed very small, resp. where the length of $\beta$
is very small, we are guaranteed that the length of $\beta$, resp.
$\alpha$, is very long.

One notes that the final part is just the result of the Baire
category theorem applied to the sets $\deform \backslash
\equal(\alpha,\beta)$. This completes the proof.
\end{proof}

In some sense the complement of $\equal$, $\Neq$, is a big set
(Baire). Another convenient way of thinking of this is that for
the generic surface in Teichm\"uller space, the length map from
the set of simple geodesics to the positive reals is injective. It
is also interesting to contrast the above lemma with the following
observation due to Randol \cite{ra80}: there are pairs of closed
geodesics $\alpha \neq \beta$ which have the same length over the
whole of Teichm\"uller space. By the above, these geodesics
necessarily have self intersections.

Lemma \ref{lem:list} proves the first part of Theorem
\ref{thm:noarc}, and the remaining part is the goal of the next
section.

\section{Total disconnectedness}\label{sec:noarc}

In this section we establish that $\Neq$, the complement of
$\equal$ in Teichm\"uller space $\deform$ is totally disconnected
i.e it contains no non-trivial arc. An {\it arc} is the continuous
image of the closed unit interval, an arc is \textit{trivial} if
it is a constant map on the interval; obviously an arc is
non-trivial if and only if it contains a subarc with distinct
endpoints. Our proof consists of applying the {\it intermediate
value theorem} to a certain continuous function which we concoct
using proposition \ref{prop:ratios}.

Let $\arc$ be an arc in Teichm\"uller space with distinct
endpoints $X,Y$. By corollary \ref{cor:ratio},
 there is a pair of simple closed geodesics
$\alpha,\beta$ such that
$$\length_X (\alpha) / \length_X (\beta) \neq \length_Y (\alpha) / \length_Y (\beta).$$
Now by proposition \ref{prop:ratios}, we see that there is an
integer $i$ such that the numbers $N(X)_i,\,N(Y)_i$, defined in
section \ref{sec:ratios}, are different. In particular, there are
simple closed geodesics $\alpha_i$ and $\beta_n$ such that
$$\length_X (\beta_n) > \length_X (\alpha_i)$$
and
$$\length_Y(\beta_n) < \length_Y (\alpha_i).$$

Now $M\mapsto \length_M(\alpha_i)-\length_M(\beta_n)$ is a
continuous function, so applying the intermediate value theorem to
$\arc$, between the points $X$ and $Y$, yields the existence of a
surface $Z\in\arc$ so that $\length_Z (\alpha_i) =
\length_Z(\beta_n)$.\\

\section{The order in lengths of simple closed geodesics
determines a point in Teichm\"uller space}

The purpose of this section is to observe that our analysis has
produced a curious fact: for closed surfaces and one-holed tori,
the orders between lengths of simple closed geodesics determine a
point in Teichm\"uller space. To see this recall that one of the
main ingredients in our proof of Theorem \ref{thm:A} is that, in
many cases, including the case of surfaces with fixed boundary
lengths and in particular closed surfaces or surfaces with cusp
boundary (Theorem \ref{thm:mono1}) and one-holed tori (lemma
\ref{lem:torus} in the Appendix), the map from Teichm\"uller space
to the marked interior simple length spectrum is projectively
injective. When this is the case, we've shown that between any two
distinct points of Teichm\"uller space there is a pair of simple
closed geodesics $\alpha$ and $\beta$ such that their order in
length is reversed. From this observation, one obtains the
following result.

\begin{theorem}\label{thm:order}
Let $\deform$ be the Teichm\"uller space of a surface with given
boundary length or of a one-holed torus with variable boundary
length. Then the map that sends a point in $\deform$ to the
relative orders in length between simple closed geodesics is
injective.
\end{theorem}

The example of two non-isometric four-holed spheres with identical
interior marked simple length spectrum at the end of Appendix
\ref{app:A} shows that one must be very careful when trying to
generalize the above theorem to arbitrary signature with variable
boundary lengths.

\section{Application to the Markoff conjecture}

Recall that $\model$ is the set of all one-holed tori with maximal
isometry group. The mapping class group of acts transitively on
connected components of $\model$ each of which is a smooth
1-dimensional submanifold of the Teichm\"uller space of one-holed
tori $\deform_{1,1}$. We now apply our analysis of multiplicities
to this subset.

\begin{theorem}\label{thm:dense}
The subset of $\model$ which does not have the Markoff isometry
property is dense.
\end{theorem}

\begin{rem} For the same reason as was given for the set
$\equal\subset\deform$ in lemma \ref{lem:list}, 
the subset consisting of those surfaces $M_t$ not having the
Markoff isometry property is a meager set in the sense of Baire.
To do this one needs to see that $\equal(\alpha,\beta)$ meets a
connected component of $\model$ in a finite number of points. This
can be done presenting $\equal(\alpha,\beta)$ and the connected
component of $\model$ as algebraic subsets of $\mathbb{R}^3$ using
the Fricke trace calculus.
\end{rem}

\begin{proof} Consider $X$ and $Y$, two distinct points in
some connected component of $\model$. Since $\model$ is a
submanifold of $\deform_{1,1}$, it follows that there is an arc
$\arc \subset \deform$ between $X$ and $Y$. By the results of
section \ref{sec:noarc}, there is a pair of geodesics
$(\alpha,\beta)$ such that $\ell_X(\alpha)<\ell_X(\beta)$ and
$\ell_Y(\beta)<\ell_Y(\alpha)$, and thus a surface $Z\in \arc$
such that $\ell_Z(\alpha)=\ell_Z(\beta)$. Since the connected
component of $\model$ is an embedded line no two surfaces $X,Y$ in
the connected component are related by an element of the mapping
class group (i.e. a change of marking). This means that if there
is a surface $X$ in the connected component such that
$\ell_X(\alpha)\neq \ell_X(\beta)$ then then there is no surface
$X'$ (in the connected component) such that an isometry of $X'$
takes $\alpha$ to $\beta$. In particular the surface $Z$
constructed above does not have the Markoff isometry property.
This establishes the claim of the theorem.
\end{proof}

\begin{rem} The surface $Z$ in the proof of the theorem
is a one-holed torus with 12 distinct simple closed geodesics of
equal length. As mentioned in the introduction, an upper bound on
the number of simple closed geodesics of equal length on a
one-holed torus would have to be at least $12$.
\end{rem}

\section{Topology of $\equal(\alpha,\beta)$}

In this section we prove Theorem \ref{thm:B} which tells us that
our hypersurfaces are connected sub-manifolds and that their
reunion is also connected. Let us begin by showing that they are
sub-manifolds.

\begin{lemma}\label{lem:subman}
$\equal(\alpha,\beta)$ is a smooth codimension 1 submanifold.
\end{lemma}

\begin{proof}
This follows from the implicit function theorem, a surface $M\in
\equal$ is a regular point for the function $M\mapsto
\ell_{M}(\alpha) - d\ell_{M}(\beta)$ if $$d\ell_{M}(\alpha) -
d\ell_{M}(\beta) \neq 0.$$ It suffices to find a single (smooth)
deformation $\phi_{t}, \phi_{0} = M$ such that
$$d\ell_{M}(\alpha).\dot{\phi_{t}} \neq
d\ell_{M}(\beta).\dot{\phi_{t}},$$ where $\dot{\phi_{t}} =
{d\phi_{t}\over dt}|_{t=0}$. In fact since $\ell_{M}(\alpha)  =
d\ell_{M}(\beta)$ it suffices to find a deformation $\phi_{t}$
such that
$${d\ell_{M}(\alpha)\over  \ell_{M}(\alpha)}.\dot{\phi_{t}} \neq
{d\ell_{M}(\beta) \over \ell_{M}(\beta) }.\dot{\phi_{t}}.$$

We explain why  Thurston's \textit{stretch maps} provide
$\phi_{t}$ satisfying this latter condition. Stretch maps were
introduced by Thurston in \cite{th98} as maps minimising the
bi-lipschitz constant between pairs of hyperbolic structures on
the same closed surface. A stretch map $\phi_{t}$ is constructed
from a maximal geodesic lamination $\lambda$; recall that a
geodesic lamination is maximal if and only every complementary
region is either an ideal triangle or a punctured monogon. One
constructs a (partial) foliation of  complementary regions by
horocyclic segments. The hyperbolic metric induces a transverse
Holder distribution on this foliation. Roughly speaking, the
stretch map consists of multiplying this transverse distribution
by $e^{t}$. This yields a hyperbolic structure such that the
length of the maximal measured sublamination of $\lambda$ is also
multiplied by $e^{t}$.
\end{proof}

To show connectedness, we will need the following lemma which
concerns spaces homeomorphic to $\reals^n$ (such as Teichm\"uller
space).

\begin{lemma}\label{lem:PB} Let $U,V$ be open connected sets of $\reals^n$ such
that $\overline{U}\cap\overline{V}=\partial U=\partial V$ and
$\overline{U}\cup V=\reals^n$. Then $\partial U=\partial V$ is
connected.
\end{lemma}

\begin{proof}
Essentially this lemma follows from the fact that $\reals^n$ has
the {\it Phragmen-Brouwer Property} \cite{brobook} which is
outlined below. A topological space is said to have the
Phragmen-Brouwer property if, for two given disjoint closed sets
$D$ and $E$ and two points $a$ and $b$, neither $D$ nor $E$
separates $a$ and $b$, then $D\cup E$ does not separate $a$ and
$b$. (A set {\it separates} two points if the two points lie on
distinct connected components of the complement of the set.)

Denote $F=\partial U=\partial V$ and suppose $F$ is not connected,
i.e., $F=F_1\cup F_2$ where $F_1$ and $F_2$ are disjoint non-empty
closed sets. Now clearly, because $\reals^n$ has the
Phragmen-Brouwer property, either the complement of $F_1$ or $F_2$
is disconnected. Suppose then that $F_1^c=A\cup B$ where $A$ and
$B$ are disjoint and open. Now $U\subset A\cup B$ is connected
thus we can suppose that $U\subset A$. Also, $F_2\subset A\cup B$
is closed and a subset of $\partial A$, so $F_2\subset A$. The set
$V\subset A\cup B$ is also connected and $\partial V \cap A
\supset F_2\neq \emptyset$ so $V\subset A$. Now
$\overline{A}=\reals^n$ and $B=\emptyset$, a contradiction.
\end{proof}

Using this lemma and stretch maps we can show the following:

\begin{lemma}\label{lem:conhyp}
The hypersurfaces $\equal(\alpha,\beta)$ are connected.
\end{lemma}

\begin{proof}
If $\alpha$ and $\beta$ are disjoint, then they can be decomposed
into a pants decomposition and because Fenchel-Nielsen parameters
determine both the points and the topology of Teichm\"uller space,
the result is obvious.

Now for intersecting $\alpha$ and $\beta$, denote by $\equal^+$,
resp. $\equal^-$, the open subsets of $\deform$ where
$\length(\alpha)>\length(\beta)$, resp. where
$\length(\alpha)<\length(\beta)$. We begin by showing that these
subsets of $\deform$ are connected. Let $M_1,M_2\in \equal^+$.
Using a stretch map, it is possible to continuously deform a
surface such that a given simple curve's length is decreased and
more so than any other simple curve. Furthermore, its length will
go to $0$. Applying a stretch map to $\beta$ on both $M_1$ and
$M_2$, we can find two finite paths in Teichm\"uller space
$c_1,c_2$ such that $c_1(0)=M_1$, $c_2(0)=M_2$ and $c_1(1)=M_1'$,
$c_2(1)=M_2'$ where $\ell_{M_1'}(\beta),\ell_{M_2'}(\beta)<
2\arcsinh 1$. Because of the property of stretch maps, these paths
remain entirely in $\equal^+$. By looking at Fenchel-Nielsen
coordinates which include the length of $\beta$, it is easy to see
that the subset of $\deform$ such that $\ell(\beta)<2\arcsinh 1$
is connected. Note that by the collar theorem and the fact that
$\alpha$ and $\beta$ intersect, any surface with
$\ell(\beta)<2\arcsinh 1$ also has $\ell(\alpha)>2\arcsinh 1$, and
thus is contained in $\equal^+$. This shows that $\equal^+$ is
connected. Using the same argument, one shows that $\equal^-$ is
connected.

Applying lemma \label{lem:PB} to $U:=\equal^+$ and $V:=\equal^-$
(and of course $\equal(\alpha,\beta)=\partial \equal^+=\partial
\equal^-$), we obtain the result.
\end{proof}

We now need to show that the reunion of the hypersurfaces is
connected which follows from the fact that $\Neq$ is totally
disconnected and the following lemma.

\begin{lemma}\label{lem:complement}
Let $N$ be a totally disconnected set of $\reals^n$ with $n>2$.
Then $\reals^n\setminus N$ is connected.
\end{lemma}

\begin{proof}
$\reals^n$ is a Cantor manifold (a result due to Urysohn for $n=3$
and to Aleksandrov for $n>3$), meaning that a closed subset that
separates has topological dimension at least $n-1$. Of course, a
totally disconnected subset $N$ is not (necessarily) closed.
However, if $N$ separates then it contains a closed subset that
separates as follows. Let $M$ be the complement of $N$ that is not
connected, thus is equal to $(A\cap M)\cup (B\cap M)$ where $A$
and $B$ are open non-empty sets which are disjoint on $M$. Thus
$A\cap B \subset N$. But as N is totally disconnected, $A\cap B =
\emptyset$. Now the complement of $A\cup B$, say $F$, is closed,
separates $\reals^n$, and is contained in $N$. But its topological
dimension cannot be greater than $0$ because $N$ has topological
dimension $0$, and thus cannot separate, a contradiction.
\end{proof}

\section{Asymptotic behaviour of $\equal(\alpha,\beta)$}\label{sec:hyptor}

The goal of this section is to prove Theorem \ref{thm:hyptor}
which gives, in the particular case of a one-holed torus with
fixed boundary length, a ``coarse" description of
$\equal(\alpha,\beta)$. We show that $\equal(\alpha,\beta)$ has
two ends and that these determine a pair of points in the Thurston
boundary of Teichm\"uller space. In fact these points are the same
pair of geodesics as in the case of flat tori outlined in the
introduction. Our point of departure is the observation that, by
the collar lemma, one expects that if $M \in \equal(\alpha,\beta)$
contains a (very) short essential simple closed geodesic $\gamma$
then
\begin{eqnarray}\label{geometric intersection}
\ii(\gamma,\alpha) = \ii(\gamma,\beta)
\end{eqnarray}

We begin by studying this equation.

\subsection{Curves and Homology}

Let $\deform^\varepsilon_{1,1}$ denote the Teichm\"uller space of
one-holed tori with some fixed boundary length $\varepsilon$. Let
$\gamma$ be an {\it oriented} primitive essential simple closed
geodesic and let $(t,\ell)$ be Fenchel-Nielsen parameters for this
choice of simple closed curve. Associated to $\gamma$ is a
foliation of $\deform^\varepsilon_{1,1}$ by level sets of
$\ell_M(\gamma)=\ell$ for $M\in \deform^\varepsilon_{1,1}$. We
make the convention that if a curve is a $k$ iterate of a
primitive curve, then each intersection point (geometric or
algebraic) is counted $k$ times.

Let us briefly recall some facts about homology classes and
oriented simple closed curves. Write $[\gamma] \in
H_{1}(M,\mathbb{Z})$ for the integer homology class determined by
a simple closed curve $\gamma$ and  $\ii_a(\cdot,\cdot)$ for the
algebraic intersection number (relative to a fixed orientation of
the surface). There is an essential simple closed curve $\gamma'$
which meets $\gamma$ just once and, after possibly changing the
orientation, we may suppose that $\ii_a(\gamma,\gamma')=1$. We say
that the pair $([\gamma],[\gamma'])$ form a\textit{ canonical
basis} for homology. The intersection numbers
$\ii_a(\gamma,\alpha)$ and $\ii_a(\delta,\alpha)$ define a unique
oriented simple closed curve up to isotopy. To each pair $(k,l)
\in \Z \times \Z$, we associate the unique oriented (not
necessarily primitive) simple closed geodesic $[\alpha]=k [\gamma]
+ l [\gamma']$ where $k=\ii_a(\alpha,\gamma')$ and
$l=-\ii_a(\alpha,\gamma)$. Recall that, in the particular case of
the one-holed torus, two curves are homologous if and only if they
are isotopic. The following equation relating algebraic and
geometric intersection numbers will of use in studying equation
(\ref{geometric intersection}):
\begin{eqnarray}\label{algebraic geometric}
|\ii_a(\alpha,\beta)|=\ii(\alpha,\beta).
\end{eqnarray}
Geometrically this means that a pair of oriented simple closed
geodesics $\alpha,\beta$ on a torus always intersect ``in the same
way", i.e., there are no arcs of $\alpha$ that leave $\beta$ and
come back to $\beta$ on the same side of $\beta$. Essentially,
equation (\ref{algebraic geometric}) holds because on a one-holed
torus no pair of arcs of simple closed geodesics form a bi-gon nor
even a bi-gon surrounding the boundary curve of the torus
\cite{busem88},\cite{ha86}.\\ 

\begin{lemma}\label{lem:2curves}
If $\cc$ is an essential simple closed geodesic such that
$\ii(\alpha,\cc)=\ii(\beta,\cc)$, then
$\pm[\cc]=[\alpha]\pm[\beta]$.
\end{lemma}

\begin{proof}
Let $[\gamma], [\gamma']$ be a basis of the integer homology, then
there exists $a,b,c,d\in \mathbb{Z}$ such that $[\alpha] = a
[\gamma] + b [\gamma']$ and $[\beta] = c [\gamma] + d [\gamma']$
with $ad-bc\neq 0$.
 Suppose $[\cc ] = x [\gamma]
+ y [\gamma']$ with $x,y\in \Z^2\setminus \{0,0\}$ then it
satisfies
$$\ii(\cc, \alpha) = \vert ay - xb \vert = \vert cy - dx \vert= \ii(\cc, \beta),$$
since geometric intersection number is just the absolute value of
the algebraic intersection number. Dropping absolute values
$$ay - bx = \pm(cy - dx)$$
that is $(x,y)$ is a solution of one of the following
\begin{eqnarray*}
(a+c)y - (b+d)x &=&0, \\
(a-c)y - (b-d)x &= & 0.\\
\end{eqnarray*}

The first equation gives $\pm[\cc]=[\alpha]+[\beta]$ and the
second $\pm[\cc]=[\alpha]-[\beta]$.
\end{proof}

\subsection{Existence}

Given a pair of primitive simple closed geodesics $\alpha$ and
$\beta$, the preceding lemma yields a pair of primitive essential
simple closed geodesics such that the geometric intersection
number with both $\alpha$ and $\beta$ is the same. We will now
show that the ends of $\equal(\alpha,\beta)$ meet the Thurston
boundary in these points Let $\gamma, \gamma'$ be a pair of simple
closed geodesic such that $([\gamma],[\gamma'])$ is a canonical
homology basis. Thus $[\alpha]=a [\gamma'] \pm n [\gamma]$ and
$[\beta]=b [\gamma'] \pm n [\gamma]$. After possibly changing
orientations of the curves $\alpha$ and $\beta$, we may suppose
that $[\alpha]=a [\gamma] - n [\gamma']$ and $[\beta]=b [\gamma] -
n [\gamma']$.

For any given torus, in \cite{mcri952} it is shown that $[\gamma]
\mapsto \ell(\gamma)$ extends to a norm, equal in fact to the
\textit{stable norm}, on the first homology group of the torus
$H_{1}(T,{\mathbb R})$. The unit ball of a norm is a convex set
and for $\gamma$ any simple closed geodesic
$[\gamma]/\ell(\gamma)$ is a boundary point of the unit ball of
the stable norm. Denote $L_\ell$  the leaf of
$\deform_{1,1}^\varepsilon$ for which  $\ell(\gamma)=\ell$.

\begin{lemma}\label{lem:hyptori}
Let $\alpha,\beta$ be a pair of distinct geodesics that meet a
simple geodesic $\gamma$ in the same number of points. Then for
any $\ell>0$, there exists a torus $T\in L_\ell$ such that
$\ell_T(\alpha) = \ell_T(\beta).$
\end{lemma}

\begin{proof}
The level sets of the function $f_{\alpha \beta}$ (defined in
lemma \ref{lem:unicity}) are just the Fenchel-Nielsen twist orbits
for $\gamma$ and thus are connected. We think of $Z$ as being a
zero of the function
$$M \mapsto \ell_{M}(\alpha) - \ell_{M}(\beta)$$

where $M$ varies in $L_\ell$.\\

To prove the existence of a solution $Z$ we shall apply the
intermediate value theorem on $L_\ell$ to our function $f_{\alpha
\beta}$. Let $D_{\gamma}$ be a positive Dehn twist along $\gamma$.
The Dehn twist acts on both the homotopy class of simple curves
and on Teichm\"ller space, sending a surface $M$ to
$D_{\gamma}.M$. For our function this means that
$$f_{\alpha \beta}(D_\gamma M) =
\ell_M (D_\gamma^{-1}(\alpha)) - \ell_M(D_\gamma^{-1}(\beta)).$$

We show that our function changes sign as $M$ varies over a twist
orbit by considering the corresponding problem for lengths of
curves on a fixed surface $M$; that is, we show that, as $k$
varies over $\mathbb Z$, there is a change of sign of
$$\ell_M (D_\gamma^{-k}(\alpha)) - \ell_M(D_\gamma^{-k}(\beta)).$$

It is convenient to suppose that $M$ is ``rectangular", that is it
admits two non-commuting reflections. As before, we choose
$\gamma'$ so that
 $([\gamma],[\gamma'])$ is a canonical basis of integer homology $H_{1}(T,{\mathbb
 Z})$
and identify $H_{1}(T,{\mathbb R})$ with ${\mathbb R}^2$ with
coordinates $(x,y)$ in the obvious way, i.e., by sending
$[\gamma],[\gamma' ] $ to the usual orthonormal basis. The two
reflections of $M$ induce automorphisms of ${\mathbb R}^2$,
reflections in respectively the $x$ and $y$ axes, and
$\bar{B}_{1}$, the closed unit ball of the stable norm is
invariant under these. Using this invariance one sees that, as we
vary over $\bar{B}_{1}$, the maximum $y$ value is attained at ${1
\over \ell_{\gamma'}}[\gamma']$. Moreover

\noindent \textbf{Observation:} \textit{ If $[\delta],[\delta'] $
on the boundary of $ \bar{B}_{1}$ are in the positive quadrant of
${\mathbb R}^2$ then the slope of the line joining these points is
negative.}

Suppose that $\alpha,\beta$ each meet $\gamma$ in exactly $n>0$
points. Then there exist integers $a,b$ such that
$$
[\alpha]  =  a [\gamma] + n[\gamma' ] $$
$$[\beta]  =  b [\gamma] + n [\gamma'].
$$
After swapping $\alpha,\beta$ we may suppose that $b>a$.

The Dehn twist along $\gamma$ induces an automorphism of
$H_{1}(T,{\mathbb Z})$, which extends to $H_{1}(T,{\mathbb R})$,
namely
$${D_\gamma}_{*} : x[\gamma] + y [\gamma'] \mapsto (x+y)[\gamma] + y [\gamma'].$$
To simplify notation for $k\in \mathbb N$ we set
$$\alpha_k :={D_\gamma}_*^k([\alpha]),\,\,
\beta_k : = {D_\beta}_*^k([\beta]).$$ The corresponding normalized
points
$${[\alpha_k]\over \ell(\alpha_k)},\,\,{[\beta_k]\over \ell(\beta_k)}$$
are in the boundary of the unit ball of stable norm
\cite{mcri952}. There exists $K\in {\mathbb Z} $ such that for all
$k>K$,
 both $[\alpha_k]$, $[\beta_k]$ are in the positive quadrant of ${\mathbb R}^2$.
The slope of the line segment joining $\alpha_k/\ell_{\alpha_k}$
to $\beta_k/\ell_{\beta_k}$ is
$${\ell(\alpha_k)\times b -  \ell(\beta_k) \times a
 \over \ell(\alpha_k) - \ell(\beta_k) },$$
which is positive if $\ell(\beta_k)\leq \ell(\alpha_k)$ and $b>a$.
This contradicts our observation above and so finishes the proof
of the lemma.
\end{proof}

\subsection{Uniqueness}

\begin{lemma}\label{lem:unicity}
Consider a leaf $L_\ell$ of $\deform^\varepsilon_{1,1}$ for some
fixed $\ell$. Denote by $M_t$ the surface of $L_\ell$ where the
twist parameter along $\gamma$ takes the value $t$. If $a>b$ then
the function $f_{\alpha
\beta}(t)=\ell_{M_t}(\alpha)-\ell_{M_t}(\beta)$ is a strictly
increasing function.
\end{lemma}

\begin{rem} Although $\alpha$ and
$\beta$ are oriented curves, this plays no role in the function
$f_{\alpha \beta}$.
\end{rem}
\begin{proof}
To prove the lemma, we shall show that the derivative of this
function is strictly positive.
 Clearly, it suffices to show that this holds whenever
$a=b+1$, the general case follows by induction.\\

We have $\ii_a(\alpha,\gamma)=\ii_a(\beta,\gamma)=n$ and
$\ii_a(\alpha,\gamma')-1=\ii_a(\beta,\gamma')=b$. Recall that a
simple closed geodesics of a one-holed torus intersect passes
through two of the three Weierstrass points of the torus. It
follows that any given any pair of simple closed geodesics
$\alpha$ and $\beta$ there is a Weierstrass point  $p\in
\alpha\cap \beta$. Consider the oriented path starting from $p$,
following $\beta$ until $\beta$ intersects $\gamma$, then
following $\gamma$ until the next intersection point with $\beta$
and so forth. The path closes up after $n$ such steps, when, after
having crossed $\gamma$ $n$ times.
 By construction this is a closed oriented piecewise smooth curve, say
$\tilde{\alpha}$, with self intersections occuring only at the
points of $\alpha\cap\beta$. Since none of the self intersections
of $\tilde{\alpha}$ are transverse it is homotopic to a simple
closed (smooth) curve. We now calculate the algebraic intersection
number $\ii(\tilde{\alpha},\cdot)$. By construction,
$\tilde{\alpha}$ is an Eulerian path for the directed graph with
vertices at $\alpha\cap\beta$ and edges the set of all
 the oriented arcs of both $\gamma$ and $\beta$ joining these points.
 Thus
$$\ii_a(\tilde{\alpha},\cdot)=\ii_a(\gamma,\cdot)+\ii_a(\beta,\cdot).$$
It follows that $\ii_a(\tilde{\alpha},\gamma)=
\ii_a(\gamma,\gamma) + \ii_a(\beta,\gamma)=0+n=n$ and
$\ii_a(\tilde{\alpha},\delta)= \ii_a(\gamma,\gamma') +
\ii_a(\beta,\gamma')=1+b=a$. This proves that $\tilde{\alpha}$ and
$\alpha$ determine the same homology class. On a hyperbolic
punctured torus there is at most one closed simple geodesic in
each homology class so, since $\tilde{\alpha}$ is homotopic to a
simple closed curve,
 $\tilde{\alpha}$ and $\alpha$ are homotopic.\\

 Lifting to the universal cover $\hyp$,
 we obtain a configuration of curves depicted in figure \ref{fig:lift}, where
the $\gamma^{(k)}$, $\beta^{(k)}$ and $p^{(k)}$ correspond to the
successive lifts of $\gamma$, $\beta$ and $p$. So, in the figure,
the orientation of $\beta$ is from top to bottom, and the
orientation of $\gamma$ from left to right. We obtain a copy of
$\alpha$ in $\hyp$ by joining the points $p^{(1)}$ to $p^{(2)}$ by
the unique oriented geodesic arc between them. Denote this
oriented geodesic arc by $\alpha'$. Note that $p$ may or may not
be a point of $\gamma$. Figure \ref{fig:lift} portrays the case
where $p$ is not a point of $\gamma$. If $p$ was a point of
$\gamma$, the lift $\gamma^{(n)}$ would pass through point
$p^{(2)}$, and the arc denoted $\beta^{1'}$ would not appear.

\begin{figure}[h]
\leavevmode \SetLabels
\L(.13*.82) $\gamma^{(1)}$\\
\L(.22*.675) $\gamma^{(2)}$\\
\L(.73*.32) $\gamma^{(n-1)}$\\
\L(.82*.16) $\gamma^{(n)}$\\
\L(.22*.9) $\beta^{(1)}$\\
\L(.31*.72) $\beta^{(2)}$\\
\L(.65*.25) $\beta^{(n)}$\\
\L(.73*.1) $\beta^{(1')}$\\
\L(.27*1.03) $p^{(1)}$\\
\L(.73*-.04) $p^{(2)}$\\
\endSetLabels
\begin{center}
\AffixLabels{\centerline{\epsfig{file =
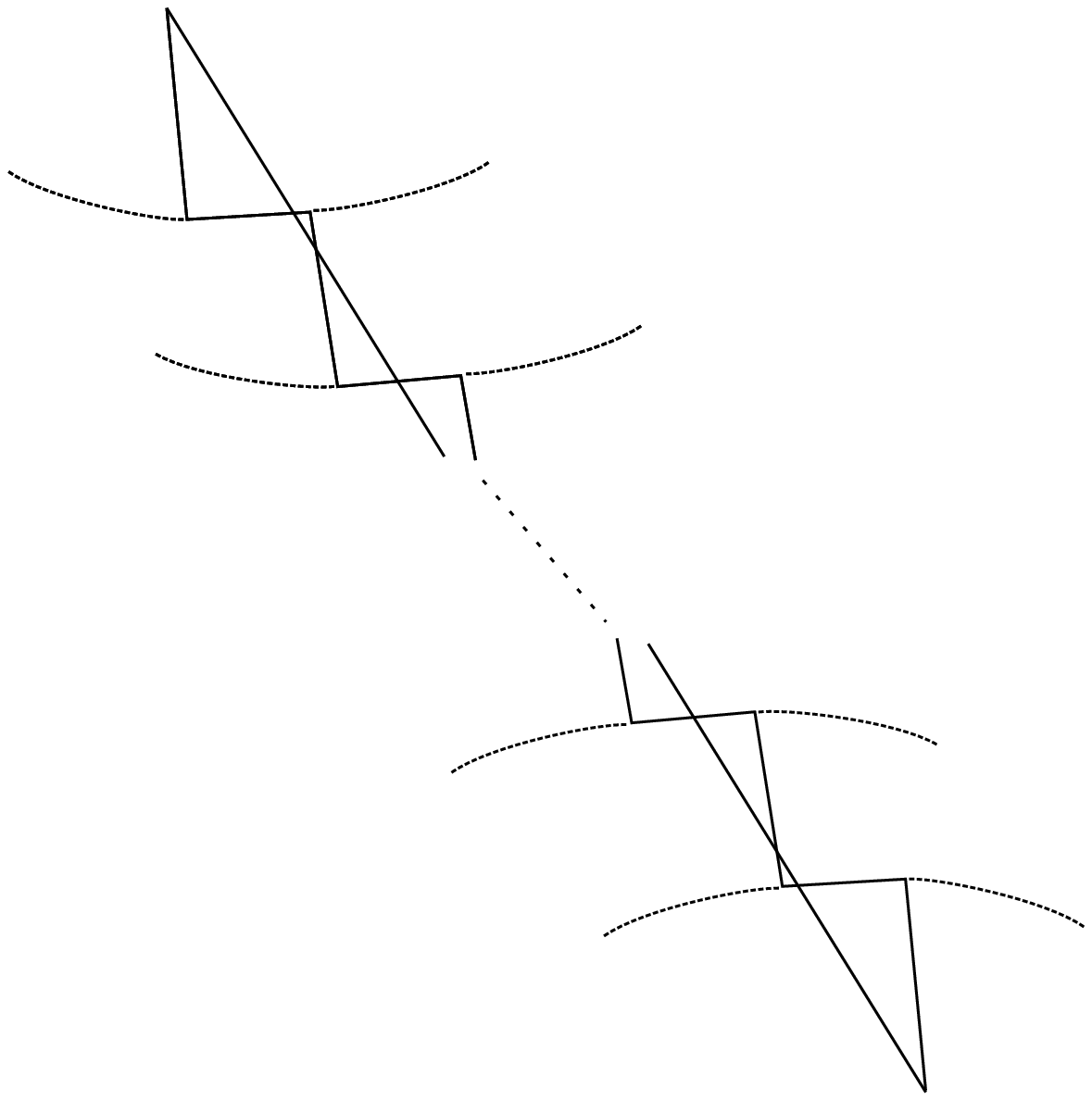,width=8cm,angle=0}}}
\end{center}
\caption{ } \label{fig:lift}
\end{figure}

We now put orders on set of  $n$ oriented angles of intersection
between $\beta$ and $\gamma$, and on the  set of  $n$ oriented
angles of intersection between $\alpha$ and $\gamma$ respectively.
Denote by $\theta_i\in ]0,\pi[$ the oriented angle between
$\beta^{(i)}$ and $\gamma^{(i)}$, and by $\tilde{\theta}_i\in
]0,\pi[$ the oriented angle between $\alpha^{(1)}$ and
$\gamma^{(i)}$. Now because $\alpha$ is homotopic to
$\tilde{\alpha}$ which is obtained by adding positive arcs along
$\gamma$ to $\beta$ we have $\theta_i> \tilde{\theta}_i$ for all
$i\in \{1,\hdots,n\}$. Now by the Kerckhoff-Wolpert formula, the
derivative along a positive twist of the function $f_{\alpha
\beta}$ is given by $\sum_{k=1}^n ( \cos \tilde{\theta}_k-\cos
\theta_k)$ which is always strictly positive as each summand $\cos
\tilde{\theta}_k-\cos \theta_k$ is strictly positive.
\end{proof}

Putting the three previous lemmas together, we have now prove
Theorem \ref{thm:C}, which rephrased in the notation introduced
above is the following.

\begin{theorem}\label{thm:hyptor}
Let $\alpha$ and $\beta$ be a pair of distinct simple closed
geodesics on a one-holed torus. The set $\equal(\alpha,\beta)$ is
a simple path joining the points in the Thurston boundary of
$\deform_{1,1}^\varepsilon$ determined by the simple closed
geodesics $[\alpha]\pm[\beta]$.
\end{theorem}

\begin{proof}
Denote by $L_\ell$ the leaf of $\deform_{1,1}^\varepsilon$ where
$\ell(\gamma)=\ell$ and by $L_\ell'$ the leaf of
$\deform_{1,1}^\varepsilon$ where $\ell(\gamma')=\ell$. By the
previous lemmas, on each leaf $L_\ell$ and $L_\ell'$, there is a
unique one-holed torus $T$ where $\ell_T(\alpha)=\ell_T(\beta)$.
As both sets $\{L_\ell\}_{\ell\in ]0,\infty[}$ foliate
$\deform_{1,1}^\varepsilon$, it follows that
$\equal(\alpha,\beta)$ is a simple path between the boundary
points of $\deform_{1,1}^\varepsilon$ $\gamma$, resp. $\gamma'$,
i.e. the point where $\ell(\gamma)=0$, respectively the point
where $\ell(\gamma')=0$.
\end{proof}

\section{Concluding remarks, problems}

In conclusion we remark that the above situation, for reasons
which should be apparent from our treatment, bears a striking
similarity to the case of the transcendental numbers. It is
relatively easy to establish the existence of aggregate of all
transcendentals as a non-empty set, and even to show that in fact
{\it most} numbers are transcendental. However, to say whether a
particular number is transcendental is a very difficult
proposition.\\

Let us end with a list of related problems.\\

\noindent{\bf Problems:}\\

\noindent\begin{enumerate}
\item It is easy to show that there exists a Euclidean torus,
with no symmetries other than the hyperelliptic involution, but
which has infinitely many distinct pairs of simple geodesics of
the same length; that is it lies on infinitely many different
hypersurfaces $\equal(.,.)$. (It is perhaps surprising that a
countable set of planes should have its set of intersections
distributed so that they clump together in such a fashion.) Is
there a hyperbolic once-punctured torus with this property?\\

\item
On the other hand is there a number $K>1$ such that if $M$ has $K$
pairs of equal length then $M$, a hyperbolic torus, must have an
isometry (other than
the hyperelliptic involution)?\\

\item
We have constructed a dense set of counterexamples to the naive
geometric generalization of the Markoff conjecture. That is, for a
dense set a of $t\geq 3$, there are matrices, in $SL(2,\reals)$,
$A_t,B_t, A_tB_t=C_t \in SL(2,\reals)$, each with trace $t$, such
that the quotient of $\hyp$ by the Fuchsian group they generate is
a torus with a hole with maximal symmmetry group but the
multiplicity of the simple spectrum is at least 12. Is it possible
to find such a $t\in \mathbb{Z}$?\\

\item As mentioned in the introduction, it is not known if simple
multiplicity is always bounded for a hyperbolic surface. This is
probably a very difficult problem. To illustrate the difficulty,
let us make the following observation. Recent results on both
upper and lower polynomial bounds on the growth of simple closed
geodesics in function of length (see \cite{mi03} and \cite{ri01})
imply that on any hyperbolic surface which is not a pair of pants,
any $\varepsilon>0$ and any $n$, there exists $n$ distinct simple
closed geodesics $\gamma_1,\hdots,\gamma_n$ such that
$|\ell(\gamma_k)-\ell(\gamma_l)|<\varepsilon$ for all
$k,l\in\{1,\hdots,n\}$. To prove this, it suffices to notice that
if this weren't true, then there would be a linear growth bound.
\end{enumerate}
\newpage
\appendix\label{app:A}
\section{Proofs of Theorems 2.1 and 2.2}

In order to prove these two theorems, we shall need some lemmas
from elementary real function theory concerning sums of
logarithmic and exponential functions. The functions we shall
study arise naturally as {\it trace polynomials} satisfied by
traces of simple closed geodesics. (Recall the trace of simple
closed geodesic $\gamma$ is the trace of a matrix in $PSL(2,\RR)$
which covers $\gamma$ and is equal to $2\cosh
\frac{\ell(\gamma)}{2}$.)

\begin{lemma}\label{lem:prelim}
Let $F:\RR^n\rightarrow \RR$ be the following function in $n>1$
variables:
$$
F(x)=\log(\sum_{k=1}^n x_k - 1)(\sum_{k=1}^n x_k - 1)-
\sum_{k=1}^n x_k \log x_k.
$$

Then $F(x)\geq 0$ for all $x$ with $\min_{k=1,\hdots,n}(x_k)\geq
1$.
\end{lemma}

\begin{proof}


We have $F(1,\hdots,1)=0$ and $\frac{\partial F(x)}{\partial
x_l}=\log (\sum_{k=1}^n x_k-1)-\log(x_l)\geq 0$ for all $l\in
\{1,\hdots,n\}$ which proves the result.
\end{proof}

\begin{lemma}\label{lem:exp}
Let $a_1,a_2,b_1,\hdots,b_n$ be positive numbers verifying
$a_1>b_k$ for all $k\in \{1,\hdots,n\}$, and let $f$ be the
following real function defined for $t\in \RR^{+}$:

$$
f(t)=a_1^t+a_2^t-\sum_{k=1}^{n}b_k^t.
$$

Then $f$ has at most one strictly positive zero.
\end{lemma}

\begin{proof}
We can suppose that $a_2=1$ as multiplying $f$ by $a_2^{-t}$ does
not change the sign of $f$, nor does it change the order of the
positive numbers.\\

Let us suppose that $f(t)\geq 0$ for some $t>0$. We shall show
that this implies $f'(t)\geq 0$. By calculation
$$
f'(t)=  a_1^t \log a_1 -\sum_{k=1}{n}b_k^t \log b_k,
$$
and thus

$$
tf'(t)= A_1 \log A_1 - \sum_{k=1}^{n} B_k \log B_k,
$$

with $A_1= a_1^t$ and $B_k=b_k^t$. Because $f(t)\geq 0$, we have
$A_1\geq \sum_{k=1}^n B_k -1$, and it follows that

$$
tf'(t)\geq \log(\sum_{k=1}^n B_k -1) (\sum_{k=1}^n B_k -1) -
\sum_{k=1}^{n} B_k \log B_k.
$$

If $B_k\geq 1$ for $k\in \{1,\hdots,n\}$, then by lemma
\ref{lem:prelim}, $tf'(t)\geq 0$. If not, then notice that

$$
tf'(t)\geq \log(\sum_{k=1}^{\tilde{n}} B_k -1)
(\sum_{k=1}^{\tilde{n}} B_k -1) - \sum_{k=1}^{\tilde{n}} B_k \log
B_k
$$

with ${\tilde{n}}<n$ and all $B_k\geq 1$ for $k\leq {\tilde{n}}$
is positive if $\tilde{n}\geq 2$. Now if $\tilde{n}\leq 1$ then
$a_1>b_1$ and $a_2=1>b_2$ and the result is trivial.
\end{proof}

\begin{corollary}\label{cor:cosh}
Let $a_1,a_2,b_1,\hdots,b_n$ be positive numbers verifying
$a_1>b_k$ for all $k\in \{1,\hdots,n\}$, and let $f$ be the
following real function defined for $t\in \RR^{+}$:

$$
f(x)=\cosh(a_1 x)+\cosh(a_2x)-\sum_{k=1}^{n}\cosh(b_kx).
$$
Then $f$ has at most one strictly positive zero.
\end{corollary}

\begin{proof}
Let us suppose that $a_1\geq a_2$. Let us remark that
$\lim_{x\rightarrow +\infty} f(x) = + \infty$. By calculation
$$
f^{(2p)}(x)=a_1^{2p}\cosh(a_1
x)+a_2^{2p}\cosh(a_2x)-\sum_{k=1}^{n}b_k^{2p}\cosh(b_k x)
$$

and

$$
f^{(2p+1)}(x)=a_1^{2p+1}\sinh(a_1 x)+a_2^{2p+1}\sinh(a_2
x)-\sum_{k=1}^{n}b_k^{2p+1}\sinh(b_k x).
$$

Because $a_1>b_k$ for all $k\in \{1,\hdots,n\}$, it follows that
there exists a $p$ such that $a_1^{2p} > \sum_{k=1}^n b_k^{2p}$
and it follows that for this $p$, we have $f^{(2p)}(x)>0$ for all
$x\geq 0$. We shall proceed backwards to $f$. Notice that for any
$q \in \NN$, $f^{(2q)}(0) = a_1^{2q}+a_2^{2q} - \sum_{k=1}^n
b_k^{2q}$ and that $f^{(2q+1)}(0)=0$. It follows that
$f^{(2p-1)}(x)>0$ for all $t>0$. For all $k$ less than $2p-1$,
$f^{(k)}(t)>0$ for all $t>0$, until there is a $p_1$ such that
$f^{(2p_1)}(0)<0$. (If such a $p_1$ does not exist then $f$ is a
strictly increasing function and the result is obvious.) Notice
that the existence of $p_1$ implies that the function defined as
$g(t):=f(t)=a_1^t+a_2^t-\sum_{k=1}^{n}b_k^t$ has a positive zero
between $2p_1$ and $2p-1$. By lemma \ref{lem:exp}, this is the
unique strictly positive zero of the function $g$. The unicity of
this zero implies that for $k=1,\hdots,p_1-1$ we have
$f^{(2k)}(0)<0$. Also, $f^{(2k+1)}(x)$ takes value $0$ for $x=0$,
is strictly decreasing for a while, and then is strictly
increasing and tends to infinity. In particular, $f'(x)$ behaves
this way. As $f(0)<0$, this completes the proof.
\end{proof}

These previous lemmas will be useful for proving our theorems for
surfaces of arbitrary signature, but are not required for proving
the results we require for one-holed tori. In the case of
one-holed tori, the map from Teichm\"uller space to the set of
lengths of interior simple closed geodesics {\it is} projectively
injective.

\begin{lemma}\label{lem:torus} There are four interior simple closed curves $\alpha$, $\beta$, $\gamma$, and $\delta$ of a
one-holed torus such that the map $\varphi:$
$M\mapsto(\length_{M}(\alpha),\length_{M}(\beta),\length_M(\gamma),\length_M(\delta))$
is projectively injective.
\end{lemma}

\begin{proof}
Let $M$ be a one-holed torus and let $\alpha$, $\beta$, $\gamma$,
and $\delta$ be the simple closed curves as in figure
\ref{fig:torusgreg}.

\begin{figure}[h]
\leavevmode \SetLabels
\L(.19*.85) $\alpha$\\
\L(.15*.72) $\beta$\\
\L(.1*.8) $\gamma$\\
\L(.28*.8) $\delta$\\
\L(.74*.19) $\alpha$\\
\L(.72*.35) $\beta$\\
\L(.64*.41) $\gamma$\\
\L(.85*.4) $\delta$\\
\L(.74*.87) $\alpha$\\
\endSetLabels
\begin{center}
\AffixLabels{\centerline{\epsfig{file =
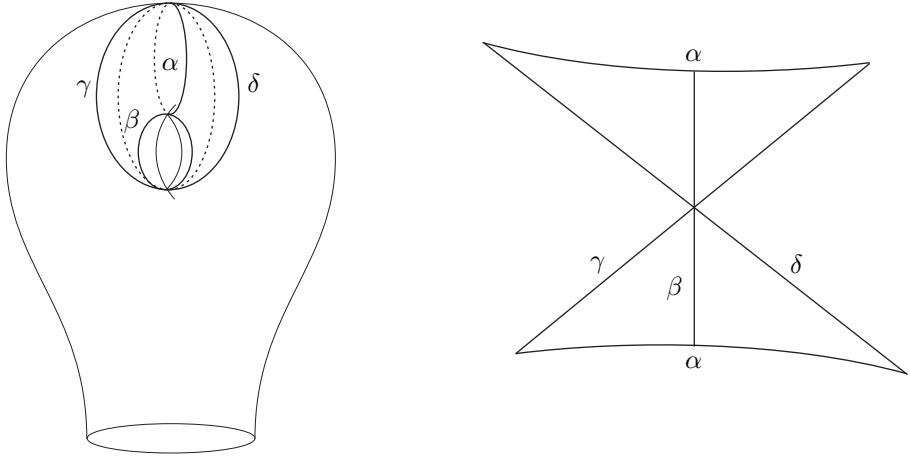,width=12cm,angle=0}}}
\end{center}
\caption{The one-holed torus with four interior geodesics and the
four curves seen in the universal cover} \label{fig:torusgreg}
\end{figure}

The remarkable fact about the geodesic representatives of simple
closed curves on a one-holed torus is that they pass through
exactly two of the three Weierstrass points of the torus in
diametrically opposite points. In the case of the curves $\alpha$,
$\beta$, $\gamma$, and $\delta$, their intersection points are all
Weierstrass points. Therefore they can be seen in the universal
cover as in figure \ref{fig:torusgreg}. It is well known that the
lengths of $\alpha$, $\beta$, and $\gamma$ determine a unique
point in Teichm\"uller space (see for example \cite{busem88}). The
lengths of these three curves up to a multiplicative constant do
{\it not} determine a unique point in Teichm\"uller space however.
For this we need the curve $\delta$. What we need to prove is that
if we have two one-holed tori $M_1$ and $M_2$ in Teichm\"uller
space with
$$
(\length_{M_1}(\alpha),\length_{M_1}(\beta),\length_{M_1}(\gamma),
\length_{M_1}(\delta))=\lambda(\length_{M_2}(\alpha),\length_{M_2}(\beta),
\length_{M_2}(\gamma),\length_{M_2}(\delta))
$$

for some $\lambda\in \RR$, then $\lambda=1$ and then, by what
precedes, $M_1=M_2$. Figure \ref{fig:torusgreg} shows four
hyperbolic triangles. Consider the two bottom ones. The side
lengths of the bottom left triangle are $\frac{\ell(\alpha)}{2}$,
$\frac{\ell(\beta)}{2}$, and $\frac{\ell(\gamma)}{2}$. The side
lengths of the bottom right triangle are $\frac{\ell(\alpha)}{2}$,
$\frac{\ell(\beta)}{2}$, and $\frac{\ell(\delta)}{2}$. The bottom
intersection point between $\alpha$ and $\beta$ forms two angles
depending on the surface $M$, say $\theta_1(M)$ and $\theta_2(M)$
such that $\theta_1 + \theta_2 = \pi$. Suppose without loss of
generality that $\lambda \geq 1$. Now if for $M_1$ the triangle
lengths are equal to $a,b,c$ and $d$, the triangle lengths for
$M_2$ are $\lambda a, \lambda b, \lambda c$ and $\lambda d$. This
implies that $\theta_1(M_1)\leq \theta_1(M_2)$ as well as
$\theta_2(M_1)\leq \theta_2(M_2)$, equality occurring only if
$\lambda =1$. As $\theta_1+\theta_2$ is always equal to $\pi$,
this concludes the proof.
\end{proof}

We are now set to prove Theorem \ref{thm:mono1}. Let us recall its
statement.

\begin{theorem}
There is a fixed finite set of simple closed geodesics
$\gamma_1,\ldots,\gamma_n$ such that the map $\varphi:$
$M\mapsto(\length_M(\gamma_1),\ldots,\length_M(\gamma_n))$ is
projectively injective.
\end{theorem}

\begin{proof}
The first step is to show that the lengths of a finite set of
simple closed geodesics determine a point in Teichm\"uller space.
Recall that a point in Teichm\"uller space is determined by the
lengths of the geodesics of a pants decomposition and twist
parameters along the interior geodesics of the pants
decomposition. In turn, the twist parameter around a geodesic is
determined by the lengths of two simple closed geodesics that
intersect the pants geodesic minimally, that mutually intersect
minimally and do not intersect any other geodesic of the pants
decomposition. (For a proof of how these lengths determine the
twist parameter, see for example \cite{bubook}). Thus the
geodesics of a pants decomposition can be completed into a finite
set of simple closed geodesics, say
$\{\gamma_1,\hdots,\gamma_m\}$, whose lengths determines a unique
point in Teichm\"uller space.\\

We shall now show that we can add an extra geodesic to one of the
finite sets described above such that the map to the lengths of
the extended set is projectively injective. There will be two
cases to consider: either the surface is of signature $(g,n)$ with
$g>0$, or not. In the first case, the surface has an embedded
one-holed torus, whose boundary curve we ensure is in the simple
closed geodesics of our pants decomposition used to create the
above set. By lemma \ref{lem:torus}, it suffices to add an simple
closed geodesic $\gamma_{m+1}$ which is interior of our torus, and
the result follows.\\

Now if $g=0$, for any pants decomposition, the surface contains an
embedded four holed sphere whose boundary geodesics are among the
pants decomposition geodesics of the surface. Denote by
$\alpha,\beta,\gamma,\delta$ these boundary geodesics, and by
$\zeta$, $\xi$ and $\upsilon$ the interior geodesics of the finite
set determined above. Consider an auxiliary geodesic $\bar{\xi}$
which is only simple closed geodesic with two intersection points
with both $\zeta$ and $\upsilon$ distinct from $\xi$. It
necessarily has four intersection points with $\xi$, and is
portrayed in figure
\ref{fig:fourholedsphere}.\\

\begin{figure}[h]
\leavevmode \SetLabels
\L(.055*.03) $\alpha$\\
\L(.04*.94) $\beta$\\
\L(.44*.03) $\gamma$\\
\L(.42*.94) $\delta$\\
\L(.17*.44) $\zeta$\\
\L(.25*.65) $\upsilon$\\
\L(.17*.79) $\xi$\\
\L(.555*.03) $\alpha$\\
\L(.54*.94) $\beta$\\
\L(.94*.03) $\gamma$\\
\L(.92*.94) $\delta$\\
\L(.67*.44) $\zeta$\\
\L(.75*.65) $\upsilon$\\
\L(.83*.71) $\bar{\xi}$\\
\endSetLabels
\begin{center}
\AffixLabels{\centerline{\epsfig{file =
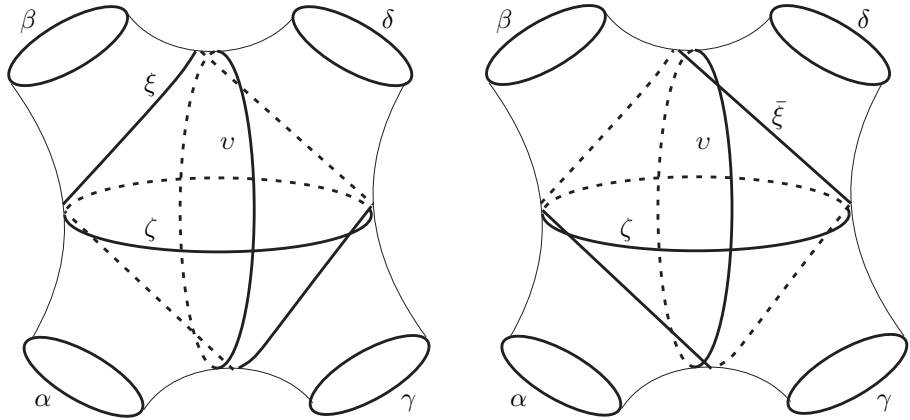,width=12cm,angle=0}}}
\end{center}
\caption{The geodesics $\zeta$, $\upsilon$, $\xi$ and $\bar{\xi}$
on the four holed sphere} \label{fig:fourholedsphere}
\end{figure}

In \cite{gase05}, it is shown that the geodesics described above
verify a certain number of trace equalities. As we require
information on the lengths of the corresponding geodesics, we have
translated these equations accordingly. One of them, that is used
to describe a fundamental domain for action of the mapping class
group of the four-holed sphere, is the following:

\begin{equation}\label{eqn:basictrace}
\begin{split}\textstyle
\cosh\frac{\ell(\xi)}{2}+\cosh \frac{\ell(\bar{\xi})}{2} +
&2\left(\cosh\frac{\ell(\alpha)}{2}\cosh\frac{\ell(\delta)}{2}
+\cosh\frac{\ell(\beta)}{2}\cosh\frac{\ell(\gamma)}{2}\right)\\
\textstyle
&=2\cosh\frac{\ell(\upsilon)}{2}\cosh\frac{\ell(\zeta)}{2}.
\end{split}
\end{equation}

By hyperbolic trigonometry this becomes

\begin{equation}\label{eqn:cosh}
\begin{split}\textstyle
\cosh\frac{\ell(\xi)}{2}+\cosh \frac{\ell(\bar{\xi})}{2}
+&\textstyle \cosh\frac{\ell(\alpha)+\ell(\delta)}{2}
+\cosh\frac{\ell(\alpha)-\ell(\delta)}{2}\\
\textstyle +\cosh\frac{\ell(\beta)+\ell(\gamma)}{2}
+\cosh\frac{\ell(\beta)-\ell(\gamma)}{2} = & \textstyle
\cosh\frac{\ell(\upsilon)+\ell(\zeta)}{2}
+\cosh\frac{\ell(\upsilon)-\ell(\zeta)}{2}.\\
\end{split}
\end{equation}

Notice that for geometric reasons,
$\cosh\frac{\ell(\upsilon)+\ell(\zeta)}{2}$ is greater than all of
the other summands. Suppose now that the map on the lengths was
not projectively injective. That would imply that there exists at
least two surfaces with their lengths verifying the above
equations. In other words, the above equations are satisfied by
both $\alpha,\beta,\gamma,\delta,\zeta,\upsilon,\xi,\bar{\xi}$ and
$t_0(\alpha,\beta,\gamma,\delta,\zeta,\upsilon,\xi,\bar{\xi})$ for
some $t_0$. Now consider the function
$f(t):=\sum_{k=1}^{2}\cosh(a_k t)-\sum_{k=1}^{6}\cosh(b_kt)$ with

\begin{eqnarray*}
a_1 &:=& \frac{\ell(\upsilon)+\ell(\zeta)}{2},\\
a_2 &:=& \frac{\ell(\upsilon)-\ell(\zeta)}{2},\\
b_1 &:=& \frac{\ell(\xi)}{2},\\
b_2 &:=& \frac{\ell(\bar{\xi})}{2},\\
b_3 &:=& \frac{\ell(\alpha)+\ell(\delta)}{2},\\
b_4 &:=& \frac{\ell(\alpha)-\ell(\delta)}{2},\\
b_5 &:=& \frac{\ell(\beta)+\ell(\gamma)}{2},\\
b_6 &:=& \frac{\ell(\beta)-\ell(\gamma)}{2}.
\end{eqnarray*}

Notice that the constants $a_1,a_2,b_1,\hdots,b_6$ satisfy the
conditions of corollary \ref{cor:cosh}. The equation $f(t)=0$ has
a positive solution given by $t=1$, and by corollary
\ref{cor:cosh}, this solution is unique, i.e., $t_0=1$, which
completes the proof.
\end{proof}

Consider a four-holed sphere. For a given choice of $\zeta$,
$\upsilon$ and $\xi$, the simple closed geodesic $\bar{\xi}$ is
uniquely determined. Similarly, consider the auxiliary curves
$\bar{\upsilon}$ and $\bar{\zeta}$ as the unique simple closed
geodesics (distinct from respectively $\upsilon$ and $\zeta$) such
that $\bar{\upsilon}$ intersects both $\xi$ and $\zeta$ twice,
respectively such that $\bar{\zeta}$ intersects both $\xi$ and
$\upsilon$ twice. As before, notice that $\upsilon$ and
$\bar{\upsilon}$ intersect four times, respectively that $\zeta$
and $\bar{\zeta}$ intersect four times. We now have a set of six
interior simple closed geodesics
$\{\xi,\bar{\xi},\upsilon,\bar{\upsilon},\zeta,\bar{\zeta}\}$
whose lengths are shown in the following lemma to determine a
finite number of points in Teichm\"uller space.

\begin{lemma}\label{lem:fourholedsphere} Let $\deform$ be the
Teichm\"uller space of the four-holed sphere. Let $M\in \deform$.
Then there are a finite many of $M_k\in \deform$ such that
$\ell_M(\gamma)=\ell_{M_k}(\gamma)$ for all $\gamma\in
\{\xi,\bar{\xi},\upsilon,\bar{\upsilon},\zeta,\bar{\zeta}\}$.
\end{lemma}

\begin{proof}
The idea is to show that the lengths of the six curves
$\{\xi,\bar{\xi},\upsilon,\bar{\upsilon},\zeta,\bar{\zeta}\}$
determine a finite number of possibilities for the lengths of the
four boundary curves $\alpha$, $\beta$, $\gamma$ and $\delta$. In
the proof of the previous theorem, we have seen that the lengths
of the geodesic representatives of $\xi, \bar{\xi},\upsilon,\zeta$
verify an interesting trace equation (equation
\ref{eqn:basictrace}) involving the boundary geodesics. The sets
$\xi,\upsilon,\bar{\upsilon},\zeta$ and
$\xi,\upsilon,\zeta,\bar{\zeta}$ verify similar equations
(equations \ref{eqn:trace2}, \ref{eqn:trace3} below). All these
equations, along with an additional trace equality which we shall
detail later (equation \ref{eqn:tracepoly}), have been proved in
\cite{gase05} and are the basis of this proof. We shall use the
following abbreviations for these equations
$a=\cosh\frac{\ell(\alpha)}{2}$ etc. (In other words, the
half-trace of a curve denoted by a Greek letter is denoted by the
lowercase corresponding roman letter.) The additional trace
equation mentioned earlier is the following:
\begin{equation}\label{eqn:tracepoly}
\begin{split} a^2+b^2+c^2&+d^2+x^2+y^2+z^2 + 4abcd -1 \\&-2xyz + 2x(ad+bc) +
2y(ab+cd) + 2z(ac+bd)=0.
\end{split}
\end{equation}

The full list of equalities that we shall use is:

\begin{eqnarray}\label{eqn:trace}
ad+bc &=& yz-\frac{x+\bar{x}}{2},\label{eqn:trace1}\\
ac+bd &=& xy-\frac{z+\bar{z}}{2},\label{eqn:trace2}\\
ab+cd &=& xz-\frac{y+\bar{y}}{2},\label{eqn:trace3}\\
a^2+b^2+c^2+d^2+4abcd &=&
1-4xyz+x\bar{x}+y\bar{y}+z\bar{z}.\label{eqn:trace4}
\end{eqnarray}

(The last equation is obtained from an obvious manipulation of
\ref{eqn:tracepoly} and \ref{eqn:trace1}, \ref{eqn:trace2},
\ref{eqn:trace3}.) Notice that the right-hand side of each
equation is determined by the lengths of
$\{\xi,\bar{\xi},\upsilon,\bar{\upsilon},\zeta,\bar{\zeta}\}$.
Denote by $f_1,\hdots,f_4$ the right hand side of each each
equation, which we shall view as given constants. It now suffices
to show that for given $f_1,\hdots,f_4 \in \RR^{+*}$, there are a
finite number of $(a,b,c,d) \in (\RR^{+*})^4$ solution to the
system of equations. There are five distinct situations to
consider.\\

The first situation is when $a=b=c=d$, which implies $f_1=f_2=f_3$
and $a=b=c=d=\sqrt{\frac{f_1}{2}}$.\\

If $a=b=c$, then

\begin{eqnarray*}
a^2+ad &=& f_1 \,\,=\,\, f_2\,\,=\,\,f_3,\\
3a^2+d^2+4a^3d &=& f_4.
\end{eqnarray*}

From these equations, we can deduce a single polynomial equation
in $a$ with a finite number of solutions. Working backwards, one
obtains a finite number of solutions
$(a=b=c,\,\,d)$ to equation \ref{eqn:trace}.\\

Now if $a=b$ and $c=d$ then
\begin{eqnarray*}
a+c &=& \sqrt{f_1+f_2},\\
ac &=& \frac{f_2}{2},
\end{eqnarray*}
which trivially implies a finite (but not unique) set of solutions
of type $(a,a,c,c)$.\\

The two remaining cases (when $a=b$, $a\neq c,d$, $c\neq d$, and
when all four variables $a,b,c,d$ are distinct) are similar in
nature. We shall give the full solution to the latter case which
is the most complicated and leave the remaining case to the
dedicated reader. Note that $a,b,c,d$ distinct implies that the
constants $f_1$, $f_2$ and $f_3$ are also distinct. By
manipulating the system of equations one obtains

\begin{align}
a &= \frac{f_2 c -f_1 d}{c^2-d^2},\label{eqn:traceII1}\\
b &= \frac{f_1 c -f_2 d}{c^2-d^2},\label{eqn:traceII2}\\
(cd-f_3)(c^2-d^2)^2
&+(f_1 c-f_2 d)(f_2 c-f_1 d)= 0,\label{eqn:traceII3}\\
\begin{split}
(f_2 c-f_1 d)^2 &+(f_1 c-f_2 d)^2+(d^2+c^2)(c^2-d^2)^2 \\
&+ 4(f_2 c-f_1 d) (f_1 c-f_2 d) c
d-f_4(c^2-d^2)^2=0.\end{split}\label{eqn:traceII4}
\end{align}

Equations \ref{eqn:traceII3} and \ref{eqn:traceII4} determine two
real planar curves. There are a finite number of pairs $(c,d)$
solution to both equations if and only if their underlying
polynomials are coprime in $\RR[c,d]$. In order to determine
whether or not they polynomials are coprime, we shall consider the
polynomials as polynomials in variable $c$, resp. in variable $d$.
Note that if the two polynomials, say $P$ and $Q$, are {\it not}
coprime in $\RR[c,d]$, then they are not coprime in either
$\RR[c]$ or $\RR[d]$. We shall then calculate the {\it resultant}
of the two polynomials in $c$, resp. in $d$, which gives a
polynomial $R_c(d)$ in variable $d$, resp. $R_d(c)$ in variable
$c$. We shall see that neither of the polynomials $R_c(d)$ and
$R_d(c)$ are identically zero. This implies that $P$ and $Q$ are
coprime in both $\RR[c]$ and $\RR[d]$, which in turn implies that
they are coprime in $\RR[c,d]$. We have

\begin{equation}
P =(cd-f_3)(c^2-d^2)^2
+(f_1 c-f_2 d)(f_2 c-f_1 d),\\
\end{equation}

\begin{equation}
\begin{split}
Q = &(f_2 c-f_1 d)^2 +(f_1 c-f_2 d)^2+(d^2+c^2)(c^2-d^2)^2\\
&+4(f_2 c-f_1 d) (f_1 c-f_2 d) c d-f_4(c^2-d^2)^2.
\end{split}
\end{equation}

The resultants are calculated by calculating the determinants of
the associated Sylvester matrixes of the two polynomials. In our
case, by computation, the resultant
$R_c(d)=\sum_{k=1}^{28}\alpha_k d^k$ is a degree $28$ polynomial
with many terms, but luckily, the leading coefficient is simple
and equal to
$$
\alpha_{28}=256(f_1-f_2)^4(f_1+f_2)^4.
$$

Similarly (and by symmetry), $R_d(c)= \sum_{k=1}^{28}\beta_k d^k$,
and
$$
\beta_{28}=256(f_1-f_2)^4(f_1+f_2)^4.
$$

As $f_1 \neq f_2$ and both $f_1$ and $f_2$ are strictly positive,
this proves the result. The remaining case can be solved using the
same method.
\end{proof}

\begin{rem}
One might hope for {\it unicity} in the previous lemma (as in the
case of the one-holed torus, lemma \ref{lem:torus}), i.e., that
the lengths of our well chosen six interior simple closed curves
determine uniquely the lengths of the boundary geodesics and thus
a unique point in Teichm\"uller space. A first remark is that, by
using the trace equalities used above, the length of the six
curves
$\{\xi,\bar{\xi},\upsilon,\bar{\upsilon},\zeta,\bar{\zeta}\}$
determine the lengths of {\it all} interior simple closed
geodesics of a given surface (see \cite{gase05} for a full proof).
In spite of this, it is easy to construct examples of pairs of
marked surfaces where the lengths of
$\{\xi,\bar{\xi},\upsilon,\bar{\upsilon},\zeta,\bar{\zeta}\}$ are
equal but the surfaces represent different points in $\deform$. To
do this consider a surface $M_1$ with marked boundary lengths
$\ell_{M_1}(\alpha)= \ell_{M_1}(\beta)=\ell_1$, and
$\ell_{M_1}(\gamma)= \ell_{M_1}(\delta)=\ell_2$. Further consider
that the marked surface $M_1$ has twist parameter equal to $0$
along $\upsilon$. Now consider a surface $M_2$ with boundary
lengths $\ell_{M_2}(\alpha)= \ell_{M_2}(\beta)=\ell_2$ and
$\ell_{M_2}(\gamma)= \ell_{M_2}(\delta)=\ell_1$, satisfying
$\ell_{M_2}(\upsilon)=\ell_{M_1}(\upsilon)$, also with twist
parameter $0$ along $\upsilon$. It is not too difficult to see
that the lengths of $\bar{\upsilon}$, $\xi$, $\bar{\xi}$, $\zeta$
and $\bar{\zeta}$ are the same for both $M_1$ and $M_2$, but yet
$M_1$ and $M_2$ are clearly distinct points in $\deform$. These
examples, however, are examples of the same surface up to
isometry. What is more surprising is that one can find examples of
surfaces that are {\it not isometric}, but have the same interior
length
spectrum.\\

Consider a surface $\tilde{M}_1$ whose marked boundary lengths are
$\ell_{\tilde{M}_1}(\alpha)= \ell_{\tilde{M}_1}(\beta)=
\ell_{\tilde{M}_1}(\gamma)=2\arccosh 2,$ and
$\ell_{\tilde{M}_1}(\delta)=2\arccosh 3$. Now set our marked
interior geodesic lengths to
$$\ell_{\tilde{M}_1}(\upsilon)=\ell_{\tilde{M}_1}(\xi)=\ell_{\tilde{M}_1}(\zeta)=q,$$

where

$$q=2\arccosh(1+(293-92\sqrt{2})^{\frac{1}{3}}+(293+92\sqrt{2})^{\frac{1}{3}}).\\$$

\smallskip
This surface has the same interior length spectrum as the surface
$M_2$ with
$$\ell_{\tilde{M}_2}(\upsilon)=\ell_{\tilde{M}_2}(\xi)=\ell_{\tilde{M}_2}(\zeta)=x,$$

\smallskip
but with $\ell_{\tilde{M}_2}(\alpha)= \ell_{\tilde{M}_2}(\beta)=
\ell_{\tilde{M}_2}(\gamma)=r,$ and $\ell_{\tilde{M}_2}(\delta)=s$
where
$$
r=\sqrt{\frac{7}{2}-\sqrt{6}},
$$

and
$$
s=\sqrt{\frac{79}{2}+15\sqrt{6}}.
$$

To see that these two surfaces have the interior simple length
spectrum, by the above it suffices to check that the left-hand
sides of equations \ref{eqn:traceII1}, \ref{eqn:traceII2}, and
\ref{eqn:traceII3} are the same, and they are by calculation. This
example is in no way isolated, and the trick to finding it was to
search for two surfaces with distinct boundary lengths, but each
with three equal boundary lengths. With this method you get at
least a real dimension $1$ family of such surfaces.
\end{rem}

We are now ready to prove Theorem \ref{thm:mono2}.

\begin{theorem}
Let $\deform$ be the Teichm\"uller space of given signature. Then,
there is a set $\Gamma:=\{\gamma_1,\hdots,\gamma_n\}$ of {\it
interior} simple closed curves such that for any given point in
$M_1\in\deform$, there are a finite number of $M_2\in\deform$ such
that $\lambda (\ell_{M_1}(\gamma_1),...,\ell_{M_1}(\gamma_n) =
(\ell_{M_2}(\gamma_1),...,\ell_{M_2}(\gamma_n))$ for some
$\lambda\in\RR$.
\end{theorem}

\begin{proof}
This is obviously true by Theorem \ref{thm:mono1} for closed
surfaces and by lemma \ref{lem:torus} for the one-holed torus. In
the remaining cases, consider the set of geodesics
$\Gamma':=\{\gamma_1,\hdots,\gamma_m\}$ constructed for the proof
of Theorem \ref{thm:mono1}. Each boundary curve (of the base
surface) in the set $\Gamma'$ is a boundary curve of a four-holed
sphere formed by elements of $\Gamma'$. To $\Gamma'$, we need to
add auxiliary curves to obtain a set $\Gamma''$ such that every
boundary curve of $\Gamma'$ is the boundary curve of a four-holed
sphere whose interior curves
$\{\xi,\bar{\xi},\upsilon,\bar{\upsilon},\zeta,\bar{\zeta}\}$ are
elements of $\Gamma''$. Now consider the set $\Gamma'''$ is
obtained by removing the boundary curves from $\Gamma''$. By the
previous lemma and the initial considerations in the proof of
Theorem \ref{thm:mono1}, the lengths of the elements of
$\Gamma'''$
determine a finite number of marked surfaces in $\deform$.\\

We now need to show that the equality $$\lambda
(\ell_{M_1}(\gamma_1),...,\ell_{M_1}(\gamma_n))=(
(\ell_{M_2}(\gamma_1),...,\ell_{M_2}(\gamma_n))$$ can only be true
for a finite number of $\lambda$ for some choice of
$\Gamma:=\{\gamma_1,\hdots,\gamma_n\}$. If the genus of the
underlying surface is not $0$, and as we have mimicked the set of
geodesics in Theorem \ref{thm:mono1}, by lemma \ref{lem:torus} it
follows that $\lambda=1$. Otherwise, let us consider one of the
four-holed spheres whose interior curves
$\{\xi,\bar{\xi},\upsilon,\bar{\upsilon},\zeta,\bar{\zeta}\}$ are
in $\Gamma'''$. Consider the further set of interior curves
obtained by a single right Dehn twist around $\zeta$ of the curves
$\xi$, $\bar{\xi}$, $\upsilon$ and $\zeta$. Notice the curve
obtained from $\zeta$ this way is of course itself, and the curve
obtained by $\xi$ is $\bar{\xi}$. Denote by $\bar{\xi}'$ and
$\upsilon'$ the images of $\bar{\xi}$ and $\upsilon$. We can now
fix $\Gamma:=\Gamma'''\cup \{\bar{\xi}',\upsilon'\}.$ Recall that
the lengths of $\xi$, $\bar{\xi}$, $\upsilon$ and $\zeta$ verify
equality \ref{eqn:cosh}. As a Dehn twist along $\zeta$ does not
change the equality satisfied by the set of curves, this gives a
second equality which, with equation \ref{eqn:cosh}, gives:
$$
\textstyle\cosh\frac{\ell(\xi)}{2}+\cosh\frac{\ell(\upsilon')+\ell(\zeta)}{2}+\cosh\frac{\ell(\upsilon')-\ell(\zeta)}{2}
$$
$$
\textstyle=\cosh\frac{\ell(\bar{\xi}')}{2}+\cosh\frac{\ell(\upsilon)+\ell(\zeta)}{2}+\cosh\frac{\ell(\upsilon)-\ell(\zeta)}{2}.
$$

\bigskip Now suppose $\ell_{M_1}(\Gamma)=\lambda
\ell_{M_2}(\Gamma)$. This implies that $f(1)=0$ and $f(\lambda)=0$
for

$$
f(t):= \sum_{k=1}^{3} \cosh(a_k t) - \sum_{k=1}^{3} \cosh(b_k t),
$$

with

\begin{eqnarray*}
a_1 &=& \frac{\ell_{M_1}(\xi)}{2},\\
a_2 &=& \frac{\ell_{M_1}(\upsilon')+\ell_{M_1}(\zeta)}{2},\\
a_3 &=& \frac{\ell_{M_1}(\upsilon')-\ell_{M_1}(\zeta)}{2},\\
b_1 &=& \frac{\ell_{M_1}(\bar{\xi}')}{2},\\
b_2 &=& \frac{\ell_{M_1}(\upsilon)+\ell_{M_1}(\zeta)}{2},\\
b_3 &=& \frac{\ell_{M_1}(\upsilon)-\ell_{M_1}(\zeta)}{2}.
\end{eqnarray*}

Unlike in corollary \ref{cor:cosh}, the function $f$ may have more
than one positive zero, but the number of zeros is clearly finite,
which implies a finite number of possible $\lambda$. This
completes the proof.
\end{proof}

\bibliographystyle{plain}
\def\cprime{$'$}

\end{document}